\documentclass[final]{siamltex}
\usepackage{helvet}
\usepackage{subeqn}
\usepackage[dvips]{epsfig,color,graphicx}
\usepackage[square]{natbib}

\newcommand{\domg}{\,d\Omega}
\newcommand{\fracd}[2]{\frac{\displaystyle{#1}}{\displaystyle{#2}}}
\newcommand{\intr}{\int_R}
\newcommand{\intK}{\int_{\raisebox{-0.7ex}{$\scriptstyle{\|{\bf k}\|\le K}$}\hspace{-4.0ex}}}
\newcommand{\into}{\int_\Omega}
\newcommand{\intinf}{\int_{-\infty}^{\infty}}

\newcommand{\DDlta}{D(\Delta)}
\newcommand{\Dlmlmp}{D_{lm,l'm'}} \newcommand{\Drhrh}{D(\rhat,\rhat)}
\newcommand{\Drhrhp}{D(\rhat,\rhat')} 
\newcommand{\grhp}{g(\rhat')} \newcommand{\Lpot}{(L+1)^2}
\newcommand{\sdsind}{\sqrt{\frac{\Delta}{\sin\Delta}}}
\newcommand{\suml}{\sum\limits}

\newcommand{\sumlm}{\sum\limits_{lm}}
\newcommand{\sumlmp}{\sum\limits_{l'm'}}
\newcommand{\sumsh}{\suml_{l=0}^{\infty}\suml_{m=-l}^{l}}
\newcommand{\sumshL}{\suml_{l=0}^{L}\suml_{m=-l}^{l}}
\newcommand{\tlofp}{\left(\frac{2l+1}{4\pi}\right)}
\newcommand{\Plm}{P_{lm}} 
\newcommand{\sff}{{\textsf{\small f}}}
\newcommand{\sfg}{{\textsf{\small g}}}
\newcommand{\sfh}{{\textsf{\small h}}}
\newcommand{\sfD}{{\textsf{\small D}}}
\newcommand{\sfG}{{\textsf{\small G}}}
\newcommand{\sfW}{{\textsf{\small W}}}
\newcommand{\Xlm}{X_{lm}}

\newcommand{\Ylm}{Y_{lm}}
\newcommand{\Ylmrh}{Y_{lm}(\rhat)}
\newcommand{\Ylmrhp}{Y_{lm}(\rhat')} \newcommand{\Ylmp}{Y_{l'm'}}

\newcommand{\rhat}{\mbf{\hat{r}}}

\newcommand{\bmp}{\begin{minipage}} \newcommand{\emp}{\end{minipage}}
\newcommand{\be}{\begin{equation}} \newcommand{\ee}{\end{equation}}
\newcommand{\bt}{\begin{tabular}} \newcommand{\et}{\end{tabular}}
\newcommand{\btstar}{\begin{tabular*}}
\newcommand{\etstar}{\end{tabular*}}
\newcommand{\btx}{\begin{tabularx}} \newcommand{\etx}{\end{tabularx}}
\newcommand{\bn}{\begin{enumerate}} \newcommand{\en}{\end{enumerate}}
\newcommand{\bd}{\begin{description}}
\newcommand{\ed}{\end{description}}
\newcommand{\bdc}{\begin{document}} \newcommand{\edc}{\end{document}}
\newcommand{\ber}{\begin{eqnarray}} \newcommand{\eer}{\end{eqnarray}}
\newcommand{\barray}{\begin{array}} \newcommand{\earray}{\end{array}}
\newcommand{\mbf}{\mathbf} \newcommand{\mrm}{\mathrm}
\newcommand{\mtl}{\mathit} \newcommand{\pl}{\partial}
\newcommand{\nb}{\nabla} \newcommand{\al}{\alpha}
\newcommand{\nnr}{\nonumber} \newcommand{\fns}{\footnotesize}
\newcommand{\nns}{\normalsize} \newcommand{\sst}{\scriptstyle}
\newcommand{\ssts}{\scriptscriptstyle} \newcommand{\rar}{\rightarrow}
\newcommand{\sct}{\section} \newcommand{\ssec}{\subsection}
\newcommand{\lbl}{\label} \newcommand{\Rar}{\Rightarrow}
\newcommand{\ssG}{\mathcal{G}} \newcommand{\ssH}{\mathcal{H}}
\newcommand{\ssL}{\mathcal{L}} \newcommand{\ssR}{\mathcal{R}}
\newcommand{\ssS}{\mathcal{S}}
\begin{document}

\title{Spatiospectral concentration on a sphere\thanks{Received by the
    editors \today.}}

\author{Frederik J. Simons\thanks{Department of Geosciences, Princeton
        University, Princeton NJ 08544, U.S.A. Now at: Department of
        Earth Sciences, University College, London WC1E 6BT,
        United Kingdom.} \and F. A. Dahlen\thanks{Department of
        Geosciences, Princeton University, Princeton NJ 08544, U.S.A.}
        \and Mark A. Wieczorek\thanks{D{\'e}partement de
        G{\'e}ophysique Spatiale et Plan{\'e}taire, Institut de
        Physique du Globe de Paris, 94701 St.~Maur-des-Foss{\'e}s},
        France.}

\label{firstpage}

\maketitle

\begin{abstract}
We pose and solve the analogue of Slepian's time-frequency
concentration problem on the surface of the unit sphere to determine
an orthogonal family of strictly bandlimited functions that are
optimally concentrated within a closed region of the sphere, or,
alternatively, of strictly spacelimited functions that are optimally
concentrated within the spherical harmonic domain. Such a basis of
simultaneously spatially and spectrally concentrated functions should
be a useful data analysis and representation tool in a variety of
geophysical and planetary applications, as well as in medical imaging,
computer science, cosmology and numerical analysis. The spherical
Slepian functions can be found either by solving an algebraic
eigenvalue problem in the spectral domain or by solving a Fredholm
integral equation in the spatial domain. The associated eigenvalues
are a measure of the spatiospectral concentration. When
the concentration region is an axisymmetric polar cap the
spatiospectral projection operator commutes with a Sturm-Liouville
operator; this enables the eigenfunctions to be computed extremely
accurately and efficiently, even when their area-bandwidth product, or
Shannon number, is large. In the asymptotic limit of a small
concentration region and a large spherical harmonic bandwidth the 
spherical concentration problem approaches its planar equivalent,
which exhibits self-similarity when the Shannon number is kept
invariant. 
\end{abstract}



\pagestyle{myheadings}
\thispagestyle{plain}
\markboth{FREDERIK~J.~SIMONS, F.~A.~DAHLEN AND
M.~A.~WIECZOREK}{SPATIOSPECTRAL CONCENTRATION ON A SPHERE}

\section{Introduction}
In a classic series of papers published in the 1960's, David Slepian
and his colleagues solved a fundamental problem in communications
engineering, namely that of optimally concentrating a signal in both
the time and frequency domains
\cite[]{Landau+60,Landau+62,Slepian64,Slepian78,Slepian+60}.  The
orthogonal family of data windows, or tapers, that arise in this
context, and their multi-dimensional extensions
\cite[]{Hanssen97,Liu+92} have been used as the basis for the
multitaper method of spectral analysis \cite[]{Percival+93,Thomson82},
and for the analysis and representation of data in a wide range of
physical, computational and biomedical disciplines (e.g., geodesy,
seismology, optics, information theory, neurology and speech
recognition). Time-frequency and time-scale concentration operators
have been studied in more general settings and a variety of one- and
multi-dimensional geometries by several authors
\cite[e.g.,][]{Daubechies88,Daubechies+88,Flandrin88,Freeden+2004,Lilly+95,Narcowich+96,Olhede+2002}.

In this paper we consider the simultaneous spatial and spectral
concentration of a real-valued function of geographical position on
the surface of the unit sphere.  The spherical multitapers that we
derive here should be useful in a number of geophysical and planetary
data analysis applications; this is our primary motivation to
undertake the present study.  In particular, we note that physical
properties, such as the thickness or elastic strength of a planetary
lithosphere, can be estimated from the cross-spectral properties of
the surface topography and associated gravitational field
\cite[]{Turcotte+81}. Such data are most commonly available as
bandlimited spherical harmonic coefficients, measured by artificial
satellites or spacecraft.  In many if not most applications, planetary
curvature prohibits the use of locally flat approximations
\cite[]{Wieczorek+98}. Thus, the determination of spatially localized
estimates of planetary properties requires spatiospectral localization
methods that go beyond those available in the plane
\cite[e.g.,][]{Simons+2003a}.  Single spherical windows or tapers have
been developed and applied in a number of recent studies
\cite[e.g.,][]{Freeden+98a,Freeden+97,Kido+2003,Simons+97a}; however,
these are neither optimally concentrated, nor as reliable as an
orthogonal family of multitapers in the extraction of robust localized
statistical information from bandlimited spherical data
\cite[]{Wieczorek+2004}.

Following an initial and extremely insightful
analysis by Gr\"{u}nbaum and his colleagues
\cite[]{Grunbaum+82}, Slepian's concentration
problem on a sphere has, to our knowledge,
been revisited quite rarely, by workers interested
in geodesy \cite[]{Albertella+99}, magnetic resonance imaging of the
human brain \cite[]{Polyakov2002} and planetary spectral analysis
\cite[]{Wieczorek+2004}. Each of these studies treats a special case.
In this paper, we pose and solve the spherical spatiospectral
concentration problem in its most general form, discuss a number
of numerical implementation methods, and analyze the flat-earth
asymptotic limit, in which the spherical concentration problem
approaches the corresponding problem on a plane.

\section{Slepian concentration problem}\label{sec:Review}
We begin with a brief review of the one-dimensional,
continuous-continuous, time-frequency concentration problem. The
results are well known so we shall just articulate them without
providing any derivations; our only objective is to provide a template
for the spherical concentration problem, which we consider in the
remainder of the paper.  We use $t$ and $\omega$ to denote time and
angular frequency, respectively, and adopt a normalization convention
in which a real time-domain signal $f(t)$ and its Fourier transform
$F(\omega)$ are related by 
\begin{equation}
\label{slepian0}
f(t)=\frac{1}{2\pi}\int_{-\infty}^{\infty}
F(\omega)e^{i\omega t}\,d\omega,\qquad 
F(\omega)=\int_{-\infty}^{\infty}
f(t)e^{-i\omega t}\,dt.
\end{equation}
The specific problem considered by Slepian \cite[]{Slepian+60}
is that of optimally concentrating a strictly bandlimited
signal $g(t)$, with a spectrum $G(\omega)$ that
vanishes for frequencies $|\omega|>W$, into a time
interval $|t|\le T$. No such bandlimited signal
$g(t)$ can be completely concentrated within
a finite interval by virtue of the Heisenberg uncertainty
principle \cite[]{Flandrin99,Messiah2000}; the optimally concentrated
signal is considered to be the one with the least energy outside of
the interval: 
\begin{equation}
\lambda=\fracd{\int_{-T}^{T}g^2(t)\,dt}
{\int_{-\infty}^{\infty}g^2(t)\,dt}=\mbox{maximum}.
\label{slepian1}
\end{equation}
Bandlimited signals $g(t)$ satisfying the variational problem~(\ref{slepian1})
have spectra $G(\omega)$ that satisfy the frequency-domain
convolutional integral eigenvalue equation 
\begin{equation}
\int_{-W}^W\fracd{\sin T(\omega-\omega')}
{\pi (\omega-\omega')}\,G(\omega')\,d\omega'
=\lambda\hspace{0.1em}G(\omega),\quad |\omega|\le W.
\label{slepian2}
\end{equation}

A closely related problem is that of concentrating the spectrum
$H(\omega)$ of a strictly timelimited function $h(t)$,
that vanishes for times $|t|> T$, into a spectral interval
$|\omega|\le W$. The appropriate measure of concentration
in this case is
\begin{equation}
\lambda=\fracd{\int_{-W}^W|H(\omega)|^2\,d\omega}
{\int_{-\infty}^{\infty}|H(\omega)|^2\,d\omega}=\mbox{maximum}.
\label{slepian3}
\end{equation}
Timelimited signals $h(t)$ whose spectra satisfy the variational
problem~(\ref{slepian3}) themselves satisfy the time-domain eigenvalue equation
\begin{equation}
\int_{-T}^{T}\fracd{\sin W(t-t')}
{\pi (t-t')}\,h(t')\,dt'=\lambda\hspace{0.1em}h(t),\quad |t|\le T.
\label{slepian4}
\end{equation}
Both equations~(\ref{slepian2}) and~(\ref{slepian4}) have the same
eigenvalues \mbox{$1>\lambda_1>\lambda_2>\cdots>0$}, with associated
time-domain eigentapers $g_1(t), g_2(t),\ldots$ and $h_1(t),
h_2(t),\ldots$ which coincide within the interval
$|t|\le T$, and associated eigenspectra
$G_1(\omega),G_2(\omega),\ldots$ and $H_1(\omega),H_2(\omega),\ldots$
which coincide within the interval $|\omega|\le W$.

A change of both the independent and dependent variables
transforms both~(\ref{slepian2}) and~(\ref{slepian4})
into the same dimensionless eigenvalue equation:
\begin{equation}
\int_{-1}^{1}\fracd{\sin TW(x-x')}
{\pi (x-x')}\,\psi(x')\,dx'=\lambda\hspace{0.1em}\psi(x),\quad |x|\le 1.
\label{slepian5}
\end{equation}
Equation~(\ref{slepian5}) shows that the eigenvalues $\lambda_1,\lambda_2,\ldots$
and suitably scaled eigenfunctions $\psi_1(x),\psi_2(x),\ldots$ depend
only upon the time-bandwidth product $TW$.  The sum of the eigenvalues
is related to this product by 
\begin{equation}
N=\sum_{\alpha=1}^{\infty}\lambda_{\alpha}=\frac{2TW}{\pi}.
\label{slepian6}
\end{equation}
Because of the characteristic step shape of the eigenvalue spectrum
\cite[]{Landau65,Slepian+65}, this so-called Shannon number
\cite[]{Percival+93} is a good estimate of the number of significant
eigenvalues, or, roughly speaking, the number of signals $f(t)$ that
can be simultaneously well concentrated into a finite time interval
$|t|\le T$ and a finite frequency interval $|\omega|\le W$. 

The integral operator acting upon $\psi$
on the left side of equation~(\ref{slepian5}) commutes with a
second-order differential operator,
\begin{equation}
\label{slepian7}
{\mathcal P}=\frac{d}{dx}(1-x^2)\frac{d}{dx}-T^2W^2x^2,
\end{equation}
which arises in the separation of
the three-dimensional scalar wave equation in prolate
spheroidal coordinates \cite[]{Slepian83}.
Because of this, it is also possible to find
the scaled eigenfunctions $\psi_1(x),\psi_2(x),\ldots$
by solving the Sturm-Liouville equation
\begin{equation}
\label{slepian7.5}
\frac{d}{dx}(1-x^2)\frac{d\psi}{dx}+(\chi-T^2W^2x^2)\psi=0,
\quad |x|\le 1,
\end{equation}
where $\chi\not=\lambda$ is the associated eigenvalue.

The bandlimited prolate spheroidal eigentapers may be chosen to be
orthonormal over the infinite time interval $|t|\le\infty$ and
orthogonal over the finite interval $|t|\le T$:
\begin{equation}
\int_{-\infty}^{\infty}g_{\alpha}g_{\beta}\,dt=\delta_{\alpha\beta}
\quad\mbox{and}\quad\int_{-T}^{T}g_{\alpha}g_{\beta}\,dt=
\lambda_{\alpha}\hspace{0.1em}\delta_{\alpha\beta}.
\label{slepian8}
\end{equation}
Almost all of the above results can be extended to the analogous
spatiospectral concentration problem for functions defined
on the surface of the unit sphere. As we shall see,
this two-dimensional problem is enriched by the arbitrary
shape of the region of spatial concentration. 

\section{Preliminaries}\label{sec:Preliminaries}
\begin{figure}[h]\centering 
\rotatebox{0}{
\includegraphics[width=0.74\textwidth]{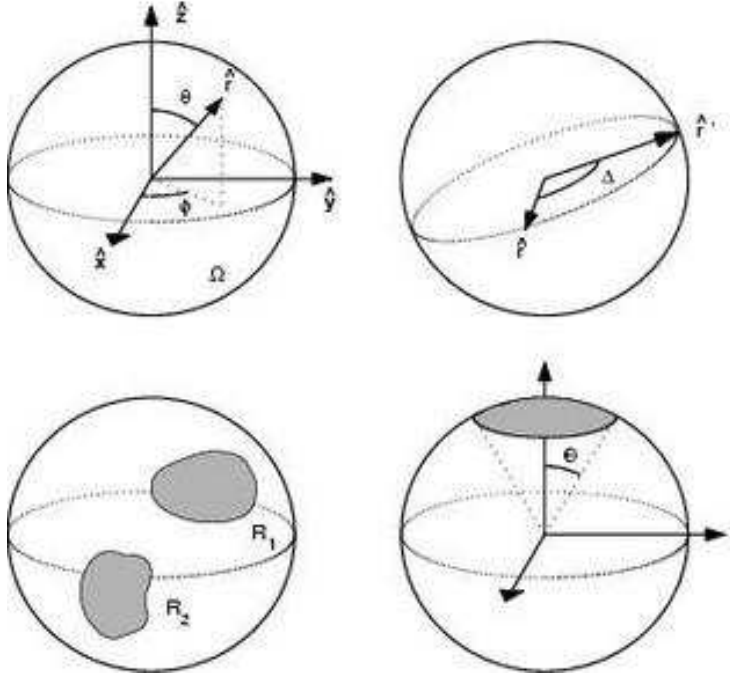}
} 
\caption{\small Sketch illustrating the geometry of the spherical
  concentration problem. Lower right shows an axisymmetric polar cap
  of colatitudinal radius $\Theta$, treated in
  Section~\ref{sec:Symmetry}. The area of the region of concentration, $R=R_1\cup
R_2\cup\cdots$, is denoted by $A$.
}
\label{sdwdiagram} 
\end{figure}

The geometry of the unit sphere $\Omega=\{\rhat: \|\rhat\|=1\}$
is depicted in Figure~\ref{sdwdiagram}.  We denote the colatitude of a
point $\rhat$ by $0\le\theta\le\pi$ and the longitude by $0\le\phi< 2\pi$,
so that $\rhat=(\theta,\phi)$ represents a geographical position
on the sphere.  The geodesic angular distance between two points
$\rhat$ and $\rhat'$ will be denoted by $\Delta$, where
\begin{equation}
\cos\Delta=\rhat\cdot\rhat'=
\cos\theta\cos\theta'+\sin\theta\sin\theta'\cos(\phi-\phi').
\end{equation}
We use $R$ to denote a region of $\Omega$ of area $A$, within which we
seek to concentrate a bandlimited function of position $\rhat$.  The
region may consist of a number of unconnected subregions, $R=R_1\cup
R_2\cup\cdots$, and it may have an irregularly shaped boundary, as
shown. The region complementary to $R$ will be denoted by $\Omega -
R$.

\subsection{Spherical harmonics}
Since we restrict attention to real-valued functions,
we use real surface spherical harmonics
$\Ylmrh=\Ylm(\theta,\phi)$ defined by
\cite[]{Dahlen+98,Edmonds96}
\ber
\Ylm(\theta,\phi)&=&\left\{
\begin{array}{l@{\quad\mbox{if}\hspace{0.6em}}l}
\rule[-2mm]{0mm}{6mm}\sqrt{2}X_{l|m|}(\theta)\cos m\phi & -l\le m<0\\
\rule[-2mm]{0mm}{6mm}X_{l0}(\theta)                     & m=0\\
\rule[-2mm]{0mm}{6mm}\sqrt{2}\Xlm(\theta)\sin m\phi & 0< m\le l,\\
\end{array}
\right.\label{Ylm}\\
\Xlm(\theta)&=&
(-1)^m\tlofp^{1/2}
\left[\frac{(l-m)!}{(l+m)!}\right]^{1/2}\!\Plm (\cos\theta),
\label{xlm}\\
\Plm (\mu)&=&
\frac{1}{2^ll!}(1-\mu^2)^{m/2}\left(\frac{d}{d\mu}\right)^{l+m}\!(\mu^2-1)^l.
\label{plm}
\eer
The quantity $0\le l\le\infty$ is known as the angular degree of
the spherical harmonic, and $-l\le m\le l$ is its angular order. The
$l\rar\infty$ asymptotic wavenumber associated with a harmonic of degree $l$ is
$\sqrt{l(l+1)}\approx l+1/2$ \cite[]{Brune64,Jeans23}.  The
function $\Plm(\mu)$ defined in (\ref{plm}) is the associated Legendre
function of integer degree $l$ and order $m$.  The spherical
harmonics $\Ylmrh$ are eigenfunctions of the Laplace-Beltrami
operator,
\begin{equation}
\nabla^2=\partial_{\theta}^2+\cot\theta\,\partial_{\theta}
+(\sin\theta)^{-2}\partial_{\phi}^2, 
\label{lapbelop}
\end{equation}
with associated
eigenvalues $-l(l+1)$. Our choice of the multiplicative
constants in equations~(\ref{Ylm}) and~(\ref{xlm})
orthonormalizes the harmonics on the unit sphere:   
\begin{equation}
\into\Ylm\Ylmp\domg=\delta_{ll'}\delta_{mm'}.
\label{normalization}
\end{equation}
The corresponding fixed-order orthogonality relations for
$X_{lm}(\theta)$ and $P_{lm}(\mu)$ are
\begin{subequations}
\ber
&&\int_0^{\pi}X_{lm}X_{l'm}\sin\theta\,d\theta
=\frac{1}{2\pi}\,\delta_{ll'},\\
\label{Xlmortho}
&&\int_{-1}^{1}P_{lm}P_{l'm}\,d\mu
=\frac{2}{2l+1}\frac{(l+m)!}{(l-m)!}\,\delta_{ll'}.
\label{Plmortho}
\eer
\end{subequations}
The integral of a Legendre polynomial
$P_l(\mu)=P_{l0}(\mu)$ over a cap $\cos\Theta\le\mu\le 1$ is
\cite[]{Byerly1893} 

\begin{equation}
\int_{\cos\Theta}^1P_l\,d\mu=\frac{1}{2l+1}\big[P_{l-1}(\cos\Theta)
-P_{l+1}(\cos\Theta)\big],,
\label{byerly}
\end{equation}
where $P_{-1}(\mu)=1$, and the product of Legendre functions at the
same argument is  
\ber
X_{lm}(\theta)X_{l'm}(\theta)&=&(-1)^m\!\!\sum_{n=|l-l'|}^{l+l'}
\displaystyle{\sqrt{\frac{(2n+1)(2l+1)(2l'+1)}{4\pi}}} \nonumber \\
&&{}\times\left(\!\begin{array}{ccc}
l & n & l' \\ 0 & 0 & 0\end{array}\!\right)
\left(\!\begin{array}{ccc}
l & n & l' \\ m & 0 & \!\!-m\end{array}\!\right)\!X_{n0}(\theta),
\label{legprod}
\eer where the arrays of indices are Wigner 3-$j$ symbols
\cite[]{Edmonds96,Messiah2000}. Of the numerous three-term recursion
relations involving the associated Legendre functions and their
derivatives, we shall make use of two in this paper, namely
\begin{subequations}
\label{needrecur}
\ber
(2l+1)\mu P_{lm}&=&(l-m+1)P_{l+1,m}+(l+m)P_{l-1,m},\\
(1-\mu^2)\frac{dP_{lm}}{d\mu}&=&(l+1)\mu P_{lm}-(l-m+1)P_{l+1,m}.
\label{grunshow5}
\eer
\end{subequations}
Finally, there are two relations involving sums of products of
Legendre functions evaluated at different
arguments that are useful in the discussion that follows.
The first is the well-known spherical harmonic addition theorem
\cite[]{Edmonds96}
\begin{equation}
\suml_{m=-l}^l\Ylmrh\Ylmrhp =
\tlofp P_l(\rhat\cdot\rhat'),
\label{additionSH}
\end{equation}
and the second is the Legendre version of the Christoffel-Darboux identity
\cite[]{Swarztrauber+2000,Szego75} 
\ber \lefteqn{
(\mu-\mu')\sum_{l=m}^L\,(2l+1)\!\left[\frac{(l-m)!}{(l+m)!}\right]\!
P_{lm}(\mu)P_{lm}(\mu')} \nonumber\hspace{2em}\\
&=&\fracd{(L-m+1)!}{(L+m)!}
\big[P_{L+1,m}(\mu)P_{Lm}(\mu')-P_{Lm}(\mu)P_{L+1,m}(\mu')\big].
\label{christdarboux}
\eer
An application of  L'H\^{o}pital's rule in equation~(\ref{christdarboux})
covers the case when $\mu=\mu'$.

\subsection{Functions on the sphere}
Let ${f}(\rhat)$ be a real-valued, square-integrable
function on the unit sphere $\Omega$. Any such function
can be expanded as a series of spherical harmonics:
\begin{equation}
\label{expansion}
{f}=\sumsh {f}_{lm}Y_{lm},\qquad
{f}_{lm}=\into {f}\hspace*{0.1em}Y_{lm}\domg.
\end{equation}
Equations~(\ref{expansion}) are the
spherical analogue of the one-dimensional Fourier
transform pair~(\ref{slepian0}).
The finite character of the unit sphere quantizes the
colatitudinal and longitudinal 
``frequencies'' $0\le l\le\infty$ and $-l\le m\le l$. We use a sans
serif $\sff$ to denote the ordered column vector of spherical harmonic
coefficients: 
\begin{equation}
\sff =\left(\begin{array}{c}
\vdots \\ f_{lm} \\ \vdots
\end{array}\right).
\end{equation}
The norm of a function $f(\rhat)$ in the spatial domain
will be denoted by 
\begin{equation}
\|{f}\|^2_{\Omega}=\into f^2\domg,
\label{spacenorm}
\end{equation}
and the norm of its spectral-domain equivalent
$\sff$ will be denoted by
\begin{equation}
\|\sff\|^2_{\infty}
=\sumsh{f}_{lm}^2.
\label{spectralnorm}
\end{equation}
Using this notation, Parseval's relation \cite[]{Percival+93}
can be written in the form
$\|{f}\|^2_{\Omega}=\|\sff\|^2_{\infty}$.
The power spectral density or variance per spherical harmonic degree $l$ and
per unit area of a function $f(\rhat)$ is defined by
\begin{equation}
\langle {f_l^2}\rangle=\frac{1}{2l+1}\suml_{m=-l}^{l}{f}_{lm}^2.
\label{specdens}
\end{equation}
We use $\delta(\rhat,\rhat')$ for the Dirac delta function
on the sphere, with the replication property
\begin{equation}
\into\delta(\rhat,\rhat')\hspace{-0.1em}f(\rhat')\domg'=f(\rhat).
\label{replication}
\end{equation}
We can write $\delta(\rhat,\rhat')$ in any of the alternative forms
\ber
\delta(\rhat,\rhat')&=&
(\sin\theta)^{-1}\delta(\theta-\theta')\delta(\phi-\phi')\nonumber \\
&=&\sumsh\Ylm(\rhat)\Ylm(\rhat')\nonumber \\
&=&\sum_{l=0}^{\infty}\tlofp\!P_l(\rhat\cdot\rhat').
\label{dirac}
\eer
The degree variance of a delta function is
constant; its power spectral density is white:
\begin{equation}
\langle\delta_l^2\rangle=\frac{1}{4\pi}\quad
\mbox{for all $0\leq l\leq\infty$}.
\end{equation}

\subsection{Bandlimited and spacelimited functions}
We are interested in two subspaces of the
space of all square-integrable functions
on the unit sphere $\Omega$. We use
\begin{equation}
\mathcal{S}_L=\big\{g\!: \langle g_l^2\rangle=0\;\,\mbox{for $L<l\le\infty$}\big\}
\label{SsubLdef}
\end{equation}
to denote the space of all bandlimited functions,
\begin{equation}
g=\sumshL g_{lm}Y_{lm},
\label{bandlg}
\end{equation}
with no power above a maximal spherical harmonic degree $L$,
and we use
\begin{equation}
\mathcal{S}_R=\big\{h\!: h=0\;\,\mbox{in $\Omega -R$}\big\}
\label{SsubRdef}
\end{equation}
to denote the space of all spacelimited functions $h(\rhat)$
that are strictly contained within a region $R$.
The space $\mathcal{S}_R$ is infinite-dimensional,
but the dimension of the space of bandlimited functions is
\begin{equation}
\dim{\mathcal{S}_L}=\suml_{l=0}^{L}(2l+1)=\Lpot,
\end{equation}
since the ordered column vector of spherical harmonic coefficients
\begin{equation}
\sfg =\left(\begin{array}{c}
g_{00} \\ \vdots \\ g_{LL}
\end{array}\right)
\end{equation}
associated with a function
$g(\rhat)$ of the form~(\ref{bandlg}) has $\Lpot$ entries. Spatial and
spectral seminorms analogous to~(\ref{spacenorm})
and~(\ref{spectralnorm}) are defined as
\begin{equation}
\|{f}\|^2_R=\intr {f}^2\domg,
\label{spacesemi}\qquad
\|\sff\|^2_L=\sumshL{f}_{lm}^2.
\end{equation}
Within the space $\mathcal{S}_R$ of spacelimited
functions, $\|h\|^2_R$ is a norm, and within
the space $\mathcal{S}_L$ of bandlimited functions,
$\|\sfg\|^2_L$ is a norm; more generally,
however, both $\|{f}\|^2_R$ and $\|\sff\|^2_L$
are seminorms.

\section{Concentration within an arbitrarily shaped region}\label{sec:Arbitrary}
The uncertainty principle \cite[]{Flandrin99,Messiah2000} stipulates
that no function can be both strictly spacelimited and strictly
bandlimited, i.e., no $f(\rhat)$ can lie in both subspaces
$\mathcal{S}_R$ and $\mathcal{S}_L$ simultaneously.  The objective of
this paper is to determine those bandlimited functions $g(\rhat)\in
\mathcal{S}_L$ that are as well contained within a spatial region $R$
as possible, and those spacelimited functions $h(\rhat)\in
\mathcal{S}_R$ whose power spectrum is as well concentrated within the
interval $0\leq l\leq L$ as possible.  As in the time-frequency case,
these two spatiospectral concentration problems will be shown to be
each other's duals. For brevity, in the discussion that follows, we
shall frequently use the abbreviations 
\begin{equation}
\sumsh=\sumlm^{\infty}\quad\mbox{and}\quad
\sumshL=\sumlm^L.
\end{equation}

\subsection{Spatial concentration of a bandlimited function}
To maximize the spatial concentration of a bandlimited
function $g(\rhat)\in{\mathcal S}_L$ within a region $R$,
we maximize the ratio of the (semi)norms:
\begin{equation}
\lambda=\fracd{\|g\|^2_R}{\|g\|^2_{\Omega}}
=\fracd{\intr g^2\domg}{\into^{}g^2\domg}=\mbox{maximum}.
\label{normratio}
\end{equation}
The two-dimensional variational problem~(\ref{normratio}) is analogous
to the one-dimensional problem~(\ref{slepian1}).  Here, as there, the
quantity $0<\lambda< 1$ is a measure of the spatial
concentration. Upon inserting the representation~(\ref{bandlg}) of
$g(\rhat)$ into (\ref{normratio}), and interchanging the order of
summation and integration, we can express $\lambda$ in the form \ber
\lambda&=& \fracd{\intr\left(\sumlm^L g_{lm}\Ylm\sumlmp^L
g_{l'm'}\Ylmp\right)\domg} {\into\left(\sumlm^L g_{lm}\Ylm\sumlmp^L
g_{l'm'}\Ylmp\right)\domg}\nnr\\ &=&\fracd{\sumlm^L
g_{lm}\sumlmp^L\left(\intr \Ylm\Ylmp\domg\right)g_{l'm'}} {\sumlm^L
g_{lm}\sumlmp^L\left(\into \Ylm\Ylmp\domg\right)g_{l'm'}}\nnr\\
&=&\fracd{\sumlm^Lg_{lm}\sumlmp^L \,\Dlmlmp \,g_{l'm'}} {\sumlm^L
g_{lm}^2},
\label{alpha}
\eer
where in the final step we have used the orthonormality relation
(\ref{normalization}), and defined the quadruply indexed quantity 
\begin{equation}
\Dlmlmp=\intr\Ylm\Ylmp\domg.
\label{Dlmlmpdef}
\end{equation}
Upon introducing the $\Lpot\times\Lpot$ matrix
\begin{equation}
\sfD = \left(\begin{array}{ccc}
D_{00,00} & \cdots & D_{00,LL} \\
\vdots & {} & \vdots \\
D_{LL,00} & \cdots & D_{LL,LL}
\label{Dmatrixdef}
\end{array}\right),
\end{equation}
with elements $D_{lm,l'm'}$, where $0\le l\le L$ and $-l\le m\le l$, we
can rewrite equation~(\ref{normratio}) as a classical matrix
variational problem \cite[]{Horn+90} in the space ${\mathcal S}_L$:
\begin{equation}
\lambda=\frac{\sfg^{\sf{\sst{T}}}\sfD \hspace{0.05em}
\sfg}{\sfg^{\sf{\sst{T}}}\sfg}
=\mbox{maximum}.
\label{Rayquot}
\end{equation}
Column vectors $\sfg$ that render the
Rayleigh quotient $\lambda$ in equation~(\ref{Rayquot}) stationary are
solutions of the $\Lpot\times\Lpot$ algebraic eigenvalue problem 
\begin{equation}
\sfD\hspace{0.05em}\sfg=\lambda\sfg.
\label{eigen1}
\end{equation}
Equation~(\ref{eigen1}) is the discrete spherical analogue of the
one-dimension spectral-domain equation~(\ref{slepian2}). The matrix
$\sfD$ is real, symmetric and positive definite,
\begin{equation}
\sfD ^{\sf{\sst{T}}}=\sfD\quad\mbox{and}\quad
\sfg^{\sf{\sst{T}}}
\sfD\hspace{0.05em}\sfg>0
\quad \mbox{for all $\sfg\ne\textsf{0}$},
\end{equation}
so the $\Lpot$ eigenvalues $\lambda$ and associated eigenvectors $\sfg$
are always real \cite[]{Horn+90}. We order the eigenvalues
$\lambda_1,\lambda_2,\ldots,\lambda_{\Lpot}$ and eigenvectors
$\sfg_1,\sfg_2,\ldots,\sfg_{\Lpot}$ so that
\begin{equation}
1>\lambda_1\ge\lambda_2\cdots\ge\lambda_{\Lpot }>0.
\label{eigorder}
\end{equation}
Each spectral-domain eigenvector ${\sfg}_{\alpha}$,
$\alpha=1,2,\ldots,\Lpot$ gives rise to an associated bandlimited
spatial eigenfunction $g_{\alpha}(\rhat)$, defined by
equation~(\ref{bandlg}).  The first inequality in~(\ref{eigorder}) is
strict because no bandlimited function can be completely
confined within a region $R$, and the last inequality is strict
because of the positive-definite character of the matrix $\sfD$.

The symmetry $\sfD^{\sf{\sst{T}}}=\sfD$ also guarantees that the
eigenvectors $\sfg_1,\sfg_2,\ldots,\sfg_{(L+1)^2}$
are mutually orthogonal \cite[]{Horn+90}. We choose them to be orthonormal:
\begin{equation}
\sfg_{\alpha}^{\sf{\sst{T}}}\sfg_{\beta}=\delta_{\alpha\beta}
\quad\mbox{and}\quad
\sfg_{\alpha}^{\sf{\sst{T}}}\sfD\hspace{0.05em}\sfg_{\beta}
=\lambda_{\alpha}\delta_{\alpha\beta}.
\label{normalize}
\end{equation}
The associated spatial eigenfunctions $g_1(\rhat),g_2(\rhat),
\ldots,g_{(L+1)^2}(\rhat)$ are in that case both orthonormal
over the whole sphere $\Omega$ and orthogonal over the region $R$:
\begin{equation}
\into g_{\alpha}g_{\beta}\domg=\delta_{\alpha\beta}\quad\mbox{and}\quad
\intr  g_{\alpha}g_{\beta}\domg=\lambda_{\alpha}\delta_{\alpha\beta}.
\label{orthog}
\end{equation}
The two spatial-domain relations in equation~(\ref{orthog}) are
equivalent to the corresponding matrix spectral
relations~(\ref{normalize}), and they are analogous to the
one-dimensional orthogonality relations~(\ref{slepian8}).  The
eigenfunction $g_1(\rhat)$ associated with the largest eigenvalue
$\lambda_1$ is the member of the space ${\mathcal S}_L$ of bandlimited
functions that is most spatially concentrated within the the region
$R$, the eigenfunction $g_2(\rhat)$ is the next best concentrated
function in ${\mathcal S}_L$ orthogonal to $g_1(\rhat)$ over both
$\Omega$ and $R$, and so on.

Written out in full using index notation, the matrix eigenvalue
equation~(\ref{eigen1}) is
\begin{equation}
\sumlmp^L \Dlmlmp g_{l'm'}=\lambda\hspace{0.1em}g_{lm}.
\label{fulleigen1}
\end{equation}
Upon multiplying equation~(\ref{fulleigen1}) by $Y_{lm}(\rhat)$
and summing over all $0\le l\le L$ and $-l\le m\le l$, the right
side yields $\lambda\hspace{0.05em}g(\rhat)$, and the left
can be manipulated as follows:
\ber
\sumlm^L\left(\,\sumlmp^L\Dlmlmp g_{l'm'}\right)\Ylmrh\nnr
&=&\sumlm^L\sumlmp^L\left(\intr \Ylmrhp
\Ylmp(\rhat')\domg'\right)g_{l'm'}\Ylmrh\nnr \nonumber \\
&=&\intr \left(\sumlm^L\Ylmrh\Ylmrhp\right)
\left(\,\sumlmp^Lg_{l'm'}\Ylmp(\rhat')\right)\domg'\nnr \nonumber \\
&=&\intr  \Drhrhp \,\grhp\domg',
\label{equiv}
\eer
where in the final step we have defined the bandlimited Dirac delta function
\begin{equation}
\Drhrhp=\sum_{l=0}^L\sum_{m=-l}^l\Ylmrh\Ylmrhp=
\sum_{l=0}^L\left(\frac{2l+1}{4\pi}\right)\!P_l(\rhat\cdot\rhat').
\label{banddelta}
\end{equation}
The above derivation demonstrates that the spectral-domain matrix
eigenvalue equation~(\ref{eigen1}) is equivalent to the spatial-domain
integral eigenvalue equation
\begin{equation}
\intr  \Drhrhp \,\grhp\domg'=
\lambda\hspace{0.05em}g(\rhat),
\quad\rhat\in\Omega.
\label{firsttimeint}
\end{equation}
Equation~(\ref{firsttimeint}) is a homogeneous Fredholm
integral equation of the second kind, with a symmetric,
separable kernel \cite[]{Kanwal71,Tricomi70}.
Upon inserting the representations~(\ref{bandlg})
and~(\ref{banddelta}) into equation~(\ref{firsttimeint}),
we recover the matrix equation~(\ref{eigen1}), so that the
spectral-domain eigenvalue problem for $\sfg$ and the spatial-domain
eigenvalue problem for a bandlimited $g(\rhat)\in{\mathcal S}_L$ are
completely equivalent. 

In summary, we can find an orthogonal family of bandlimited
eigenfunctions that are optimally concentrated within a region $R$ on
the unit sphere $\Omega$ either by solving the $\Lpot\times\Lpot$
matrix eigenvalue problem~(\ref{eigen1}) for the spectral-domain
eigenvectors $\sfg_1,\sfg_2,\ldots,\sfg_{\Lpot}$,
or by solving the Fredholm integral equation~(\ref{firsttimeint})
for the associated spatial-domain eigenfunctions $g_1,g_2,\ldots,g_{(L+1)^2}$.
Either method determines the optimally concentrated eigenfunctions
at all points $\rhat\in\Omega$, i.e., both in the region $R$, where
they are concentrated, and in the complementary region $\Omega -R$,
where they exhibit inevitable leakage.

\subsection{Spectral concentration of a spacelimited function}
Instead of seeking to concentrate a bandlimited function
$g(\rhat)\in{\mathcal S}_L$ within a spatial region $R$,
we may seek to concentrate a spacelimited function
$h(\rhat)\in{\mathcal S}_R$ within a spectral interval
$0\le l\le L$. A suitable measure of concentration
is then the spectral norm ratio, analogous to the one-dimensional
ratio~(\ref{slepian3}):
\begin{equation}
{\lambda}=\frac{\displaystyle{\|{\sfh}\|^2_L}}
{\|\sfh\|^2_{\ssts\infty}}=
\fracd{\sum_{l=0}^{L}\sum_{m=-l}^lh^2_{lm}}
{\sum_{l=0}^{\infty}\sum_{m=-l}^l h^2_{lm}}=\mbox{maximum}.
\label{normratio2}
\end{equation}
Upon inserting the representation of the spherical
harmonic expansion coefficients,
\begin{equation}
h_{lm}=\intr h\,Y_{lm}\,d\Omega,
\label{hexpansion}
\end{equation}
and interchanging the order of summation and integration, we can rewrite
the ratio~(\ref{normratio2}) in the form
\ber
\lambda&=&\fracd{\sumlm^L\intr h(\rhat)\Ylmrh\domg
\intr h(\rhat')\Ylmrhp \domg'}
{\sumlm^{\infty}\intr h(\rhat)\Ylmrh\domg
\intr h(\rhat')\Ylmrhp \domg'}\nnr\\
&=&\fracd{\intr\intr h(\rhat)\left(\sumlm^L
\Ylmrh\Ylmrhp \right)h(\rhat')\domg'\domg}
{\intr\intr h(\rhat)\left(\sumlm^{\infty}
\Ylmrh\Ylmrhp \right)h(\rhat')\domg'\domg}\nnr\\
&=&\fracd{\intr\intr h(\rhat)\Drhrhp h(\rhat')\domg\domg'}
{\intr h^2(\rhat) \domg},
\label{normratio3}
\eer where in the final step we have made use of replication property
(\ref{replication}) of the delta function (\ref{dirac}) and the
definition~(\ref{banddelta}) of the kernel $\Drhrhp$. Functions
$h(\rhat)\in{\mathcal S}_R$ that render the Rayleigh
quotient~(\ref{normratio3}) stationary are solutions of the Fredholm
integral eigenvalue equation
\begin{equation}
\intr\Drhrhp h(\rhat')\domg'
=\lambda\hspace{0.05em}h(\rhat),\quad\rhat\in R.
\label{eigen2}
\end{equation}
Equation~(\ref{eigen2}) is the spherical analogue of the
one-dimensional time-domain eigenvalue equation~(\ref{slepian4}). In
fact, this equation for $h(\rhat)\in{\mathcal S}_R$ is identical to
equation~(\ref{firsttimeint}) for $g(\rhat)\in{\mathcal S}_L$. The
only difference is that equation~(\ref{firsttimeint}) is applicable on
the entire sphere $\Omega$, whereas the domain of
equation~(\ref{eigen2}) is limited to the region $R$, within which
$h(\rhat)\not= 0$.  Evidently, the eigenfunctions $h(\rhat)$ that
maximize the spectral norm ratio~(\ref{normratio2}) are identical,
within the region $R$, to the eigenfunctions $g(\rhat)$ that maximize
the spatial norm ratio~(\ref{normratio}):
\begin{equation}
h(\rhat)=\left\{\begin{array}{ll}
g(\rhat) & \mbox{if $\rhat\in R$}\\
0 & \mbox{otherwise}.
\end{array}\right.
\label{hequalsg}
\end{equation}
Every one of the $(L+1)^2$ bandlimited eigenfunctions
$g_\alpha\in{\mathcal S}_L$ gives rise to a spacelimited eigenfunction
$h_\alpha\in{\mathcal S}_R$, defined by the
restriction~(\ref{hequalsg}). The associated eigenvalues
$\lambda_\alpha$ are a measure of both the spatial concentration of
the bandlimited eigenfunctions within the region $R$ and the spectral
concentration of the spacelimited eigenfunctions within the interval
$0\le l\le L$.  The fractional spatial energy $1-\lambda_{\alpha}$
leaked by an eigenfunction $g_{\alpha}(\rhat)$ to the region $\Omega -
R$ is identical to the fractional spectral energy leaked by
$h_{\alpha}(\rhat)$ into the higher degrees $L< l\le\infty$. If we had
started with the variational prescription~(\ref{normratio2}) rather
than~(\ref{normratio}), we could have obtained the integral
equation~(\ref{firsttimeint}) governing a bandlimited function
$g(\rhat)\in{\mathcal S}_L$ by simply extending the domain of solution
of equation~(\ref{eigen2}) to the whole sphere $\Omega$.

The spacelimited eigenfunctions defined by equation~(\ref{hequalsg}) are
orthogonal but not orthonormal over the whole sphere $\Omega$ and over
the region $R$: 
\begin{equation}
\label{horthog}
\into h_{\alpha}h_{\beta}\domg=\intr h_{\alpha}h_{\beta}\domg=
\lambda_{\alpha}\delta_{\alpha\beta}.
\end{equation}
The spherical harmonic coefficients $h_{lm}$, with $0\leq l\leq\infty$
and $-l\leq m\leq l$, are given in terms of those of $g_{lm}$, with
$0\leq l\leq L$ and $-l\leq m\leq l$, by
\begin{equation}
\label{hlmglm}
h_{lm}=\sum_{l'=0}^L\sum_{m'=-l'}^{l'}D_{lm,l'm'}g_{l'm'},
\end{equation}
which reduces to $h_{lm}=\lambda g_{lm}$ for $0\leq l \leq L$.  In
addition to the $\Lpot$ eigenfunctions with associated non-zero
eigenvalues $\lambda_1,\lambda_2,\ldots,\lambda_{\Lpot}$,
equation~(\ref{eigen2}) has an infinite-dimensional null space of
eigenfunctions with associated eigenvalue $\lambda=0$. Every function
$h(\rhat)$ that vanishes in $\Omega- R$ and has no energy whatsoever
in the interval $0\le l\le L$ is a member of this null space.

\subsection{Significant and insignificant eigenvalues}
\label{Section:sigeivs}
The sum of the eigenvalues of the matrix
${\textsf D}$ defined in equation~(\ref{Dmatrixdef}) is
\begin{equation}
N=\sum_{\alpha =1}^{\Lpot}\lambda_{\alpha}=
\mbox{tr}\,{\textsf D}=\sum_{l=0}^L\sum_{m=-l}^l D_{lm,lm}.
\label{tracedef}\label{Nsimple}
\end{equation}

Upon substituting for the diagonal matrix elements $D_{lm,lm}$
from equation~(\ref{Dlmlmpdef}) and making use of the spherical
harmonic addition theorem~(\ref{additionSH}), this simplifies to
\ber
N&=&\sum_{l=0}^L\sum_{m=-l}^l \intr\Ylm\Ylm\domg\nnr\\
&=&\sum_{l=0}^L\intr\left(\sum_{m=-l}^l\Ylm\Ylm\right)d\Omega\nnr\\
&=&\sum_{l=0}^L\left(\frac{2l+1}{4\pi}\right)\intr P_l(1)\domg\nnr\\
&=&\Lpot\,\frac{A}{4\pi},
\label{Nsimple2}
\eer
where in the final step we have used the identity $P_l(1)=1$
to express the result in terms of the surface area, $A=\intr\domg$, of
the region of concentration $R$.

We can alternatively obtain the result~(\ref{Nsimple}) by starting
with the spatial eigenvalue equation~(\ref{firsttimeint}).
The kernel $\Drhrhp$ can be expressed in terms of the
spatial eigenfunctions $g_1,g_2,\ldots,g_{(L+1)^2}$, which constitute
a basis for ${\mathcal S}_L$, in the form
\begin{equation}
\Drhrhp=\sum_{\alpha=1}^{\Lpot}g_{\alpha}(\rhat)g_{\alpha}(\rhat').
\label{newbasis}
\end{equation}
The representation~(\ref{newbasis}) is the spherical analogue of
Mercer's theorem \cite[]{Kanwal71,Tricomi70}. Algebraically, we can
regard the relation~(\ref{newbasis}) as having been obtained from the
original representation~(\ref{banddelta}) by an orthogonal
transformation from one basis, $\Ylm$, $0\le l\le L$ and $-l\le m\le
l$, to another, $g_{\alpha}$, $\alpha=1,2,\ldots,\Lpot$. Setting
$\rhat'=\rhat$ in equation~(\ref{newbasis}), and integrating over the
region $R$, we obtain
\begin{equation}
\intr  \Drhrh\domg=
\sum_{\alpha=1}^{\Lpot}\intr g_{\alpha}^2(\rhat)\domg=
\sum_{\alpha=1}^{\Lpot}\lambda_{\alpha}.
\label{sumlambda}
\end{equation}
Alternatively, setting $\rhat'=\rhat$ in equation~(\ref{banddelta})
and integrating over $R$ yields
\begin{equation}
\intr\Drhrh\domg=\suml_{l=0}^{L}\tlofp\intr 
P_l(1)\domg=\Lpot\,\frac{A}{4\pi}.
\label{sumlambda2}
\end{equation}
Comparing equations~(\ref{Nsimple})--(\ref{Nsimple2})
and~(\ref{sumlambda})--(\ref{sumlambda2}),
we see that we can write the sum of the eigenvalues~(\ref{tracedef})
in any of the equivalent ways
\begin{equation}
N=\sum_{\alpha=1}^{\Lpot}\lambda_{\alpha}=
\mbox{tr}\,{\textsf D}=\intr\Drhrh\domg=\Lpot\,\frac{A}{4\pi}.
\label{shannon}
\end{equation}

The quantity $N$ is the spherical analogue of the Shannon
number~(\ref{slepian6}) in Slepian's one-dimensional concentration problem.
Eigenfunctions $g_{\alpha}(\rhat)$ that are well
concentrated within the region $R$ will have associated eigenvalues
$\lambda_{\alpha}$ that are near unity, whereas those that are poorly
concentrated will have associated eigenvalues $\lambda_{\alpha}$ that
are near zero.  If, as in the one-dimensional problem,
the spectrum of eigenvalues $\lambda_1,\lambda_2,
\ldots,\lambda_{\Lpot}$ has a narrow transition band from values near
unity to values near zero, then the total number of
significant ($\lambda_{\alpha}\approx 1$) eigenvalues will be well
approximated by the (rounded) sum~(\ref{tracedef}). For this reason, we expect
$N$ to be a good estimate of the number of significant
eigenvalues. Roughly speaking, the spherical Shannon
number~(\ref{shannon}) is the dimension of the space of
two-dimensional functions $f(\rhat)$ that are both approximately
limited in the spectral domain to spherical harmonic degrees $0\le
l\le L$ and approximately limited in the spatial domain to an
arbitrarily shaped region $R$ of area $A$ \cite[]{Landau65,Landau67}.

Rather than seeking a bandlimited function $g(\rhat)\in{\mathcal S}_L$
that is optimally concentrated within a spatial region $R$, we could
have decided to seek one that is optimally excluded from $R$, i.e.,
one that is optimally concentrated within the complementary region
$\Omega -R$. In that case we would have sought to minimize rather than
maximize the Rayleigh quotient~(\ref{normratio}). In fact, all that we
have found are the bandlimited functions $g(\rhat)\in{\mathcal S}_L$
that render the eigenvalue $\lambda$ stationary, so we have actually
solved the containment problem~(\ref{normratio}) and the exclusion
problem simultaneously.  The optimally concentrated
eigenfunctions and the optimally excluded eigenfunctions are
identical, but with the ordering indices reversed, i.e., the
bandlimited function that is most excluded from $R$ and most
concentrated within $\Omega -R$ is $g_{\Lpot}(\rhat)$, with the
smallest associated eigenvalue $\lambda_{\Lpot}$.  Whenever the area
$A$ of the region $R$ is a small fraction of the area of the sphere,
$A\ll 4\pi$, there will be many more well excluded eigenfunctions with
insignificant ($\lambda_{\alpha}\approx 0$) eigenvalues, than well
concentrated eigenfunctions with significant ($\lambda_{\alpha}\approx
1$) eigenvalues, i.e., $N\ll\Lpot$.

The sum of the squares of the $\Lpot$ bandlimited eigenfunctions
$g_{\alpha}(\rhat)$ is a constant, independent
of position $\rhat$ on the sphere $\Omega$. This is another
consequence of Mercer's theorem~(\ref{newbasis}) and the definition
(\ref{banddelta}): 
\begin{equation}
\sum_{\alpha=1}^{\Lpot}g_{\alpha}^2(\rhat)
=D(\rhat,\rhat)=\frac{\Lpot}{4\pi}=\frac{N}{A}.
\label{sumofsq}
\end{equation}
If the first $N$ eigenfunctions $g_1,g_2,\ldots,g_N$ have eigenvalues
near unity and lie mostly within $R$, and the remainder
$g_{N+1},g_{N+2},\ldots,g_{(L+1)^2}$ have eigenvalues near zero and
lie mostly in $\Omega-R$, then we expect the eigenvalue-weighted sum of
squares to be
\begin{equation}
\sum_{\alpha=1}^{\Lpot}\lambda_{\alpha}\hspace{0.1em}
g_{\alpha}^2(\rhat)\approx
\sum_{\alpha=1}^N\lambda_{\alpha}\hspace{0.1em}
g_{\alpha}^2(\rhat)\approx
\left\{\begin{array}{ll}
N/A & \mbox{if $\rhat\in R$}\\
0 & \mbox{otherwise.}
\end{array}\right.
\label{sumofsq2}
\end{equation}
The terms with $N+1\le\alpha\le\Lpot$
should be negligible, so it is immaterial
whether they are included in the sum~(\ref{sumofsq2}) or not.
Taken together, the first $N$ orthogonal eigenfunctions
$g_{\alpha},\alpha=1,2,\ldots,N$,
with significant eigenvalues $\lambda_{\alpha}\approx 1$,
provide an essentially uniform coverage of the
\mbox{region $R$}. This is really the essence of
the spatiospectral concentration problem; the
number of degrees of freedom is reduced from
$\mbox{dim}\,{\mathcal S}_L=\Lpot$ to the Shannon
number $N=\Lpot A/(4\pi)$.

\subsection{Abstract operator formulation}
\label{subsec:AbsOps}
We conclude this section on the concentration problem for an
arbitrarily shaped region by reiterating the above results using an
abstract operator notation.  We use $\ssH$ to denote the operator that
acts upon square-integrable functions $f(\rhat)$ in the spatial domain
to produce the associated infinite-dimensional column vectors $\sff$
of spherical harmonic coefficients $f_{lm}$ in the spectral domain,
and we use $\ssH^{-1}$ to denote its inverse, so that
\begin{equation}
\ssH f=\sff\quad\mbox{and}\quad\ssH^{-1}\sff=f.
\label{oper1}
\end{equation}
We also introduce two projection operators, $\ssR$ and $\ssL$,
which project onto the space $\ssS_R$ of spacelimited functions
and the space $\ssS_L$ of bandlimited functions, respectively.
The first of these acts to spatially restrict functions
$f(\rhat)$ in the spatial domain,
\begin{equation}
\ssR f(\rhat)=\left\{\begin{array}{ll}
f(\rhat) & \mbox{if $\rhat\in R$}\\
0 & \mbox{otherwise,}
\end{array}\right.
\label{oper2}
\end{equation}
whereas the second acts to bandlimit column
vectors in the spectral domain,
\begin{equation}
\ssL\,\sff=\ssL\left(\begin{array}{c}
f_{00} \\ \vdots \\ f_{\infty\infty} \end{array}\right)=
\left(\begin{array}{c}
f_{00} \\ \vdots \\ f_{LL} \end{array}\right).
\label{oper3}
\end{equation}
Finally, we introduce a notation for the standard inner product
in both domains:
\begin{equation}
\langle f, f'\rangle_{\Omega}=\int_{\Omega}f\hspace{-0.1em}f'\domg
\quad\mbox{and}\quad\langle \sff, \sff^{\,\prime}
\rangle_{\infty}=\sff^{\hspace*{0.1em}\sf{\sst{T}}}\sff^{\,\prime}.
\label{oper4}
\end{equation}
Parseval's relation can be written using this notation in the form
$\langle f, f'\rangle_{\Omega}=\langle \sff, \sff^{\,\prime}
\rangle_{\infty}$.  The spatial and spectral norms introduced
in equations~(\ref{spacenorm}) and~(\ref{spectralnorm}) are
given by $\|f\|_{\Omega}^2=\langle f, f\rangle_{\Omega}$
and $\|\sff\|_{\infty}^2=\langle \sff, \sff
\,\rangle_{\infty}$.

The spatial-concentration variational problem~(\ref{normratio}) and
the spectral-concentration variational problem~(\ref{normratio2})
can be written using this operator notation in the form
\begin{subequations}
\ber
\lambda&=&\fracd{\langle\ssR\ssH^{-1}\ssL\,\sff,
\ssR\ssH^{-1}\ssL\,\sff\,\rangle_{\Omega}}
{\langle\ssH^{-1}\ssL\,\sff,
\ssH^{-1}\ssL\,\sff\,\rangle_{\Omega}}=\mbox{maximum},\\
\label{oper5}
\lambda&=&\fracd{\langle\ssL\ssH\ssR f, \ssL\ssH\ssR f\rangle_{\infty}}
{\langle\ssH\ssR f, \ssH\ssR f\rangle_{\infty}}=\mbox{maximum}.
\label{oper6}
\eer
\end{subequations}
The associated spectral and spatial-domain eigenvalue equations are
\begin{subequations}
\label{gruncomp}
\ber
(\ssL\ssH\ssR\ssH^{-1\!}\ssL)(\ssL\,\sff\,)
&=&\lambda\hspace{0.1em}(\ssL\,\sff\,),\label{oper7}\\
(\ssR\ssH^{-1\!}\ssL\ssH\ssR)(\ssR f)&=&\lambda\hspace{0.1em}(\ssR f)\label{oper8},
\eer
\end{subequations}
where we have made use of the facts that $\ssH$ and $\ssH^{-1}$ are each others'
transposes, that both $\ssR$ and $\ssL$ are their own transposes,
and that $\ssR^2=\ssR$ and $\ssL^2=\ssL$. Equations~(\ref{gruncomp})
are the operator equivalents of the algebraic
eigenvalue equation~(\ref{eigen1}) and the
integral eigenvalue equation~(\ref{eigen2}).  Any solution of
equation~(\ref{oper7}) is a bandlimited column vector
of the form $\sfg=\ssL\,\sff$, whereas
any solution of equation~(\ref{oper8}) is a spacelimited
function of the form $h=\ssR f$. Both the spectral-domain operator
$\ssL\ssH\ssR\ssH^{-1\!}\ssL$ and the spatial-domain operator
$\ssR\ssH^{-1\!}\ssL\ssH\ssR$ are symmetric by inspection.
Application of the operator product $\ssH^{-1}\ssL\ssH$ acts
to bandlimit an arbitrary function $f$, an operation referred to as
spherical or uniform-resolution filtering, or as triangular truncation,
in numerical analysis
\cite[]{Boehme+2003a,Jakob-Chien+97,OuldKaber96,Swarztrauber+2000}.

\section{Concentration within an axisymmetric polar cap}
\label{sec:Symmetry}
We turn our attention next to the special but important case in
which the region of concentration is a circularly symmetric
cap of colatitudinal radius $\Theta$, centered on the north pole,
as illustrated in the lower right of Figure~\ref{sdwdiagram}:
\begin{equation}
R=\big\{\theta: 0\le\theta\le\Theta\big\}.
\label{polarcap}
\end{equation}
In practical applications, the eigenfunctions that are optimally concentrated
within such a polar cap can be rotated to an arbitrarily positioned
circular cap on the unit sphere using standard spherical harmonic
rotation formulae \cite[]{Blanco+97,Dahlen+98,Edmonds96,Masters+98}.

\subsection{Decomposition of the spectral-domain eigenvalue problem}
The matrix elements~(\ref{Dlmlmpdef}) reduce in the case~(\ref{polarcap}) to
\begin{equation}
\Dlmlmp =
2\pi\,\delta_{mm'}\int_{0}^{\Theta}
X_{lm}X_{l'm'}\sin\theta\,d\theta.
\label{blockdia}
\end{equation}
The Kronecker delta $\delta_{mm'}$ renders the matrix $\sfD$ in
equation~(\ref{Dmatrixdef}) block-diagonal: 
\begin{equation}
\sfD = \left(\begin{array}{cccccc}
\sfD_0 &&&&& \\
& \sfD_{\pm 1} &&&&\\
&& \sfD_{\pm 1} &&&\\
&&&\ddots &&\\
&&&& \sfD_{\pm L} &\\
&&&&& \sfD_{\pm L} \end{array}\right),
\label{blockdia2}
\end{equation}
where every $(L-m+1)\times (L-m+1)$ submatrix
$\sfD_{\pm m}\not=\sfD_0$ occurs twice as a result
of the doublet degeneracy associated with $\pm m$.
Rather than solving the full
$\Lpot\times\Lpot$ eigenvalue equation~(\ref{eigen1}),
we may solve a series of smaller eigenvalue equations,
$\sfD_{\pm m}\sfg_{\pm m}=\lambda_{\pm m}\sfg_{\pm m}$,
one for each order $\pm m$.

\begin{figure}[b]\centering 
\rotatebox{0}{
\includegraphics[width=1\textwidth]{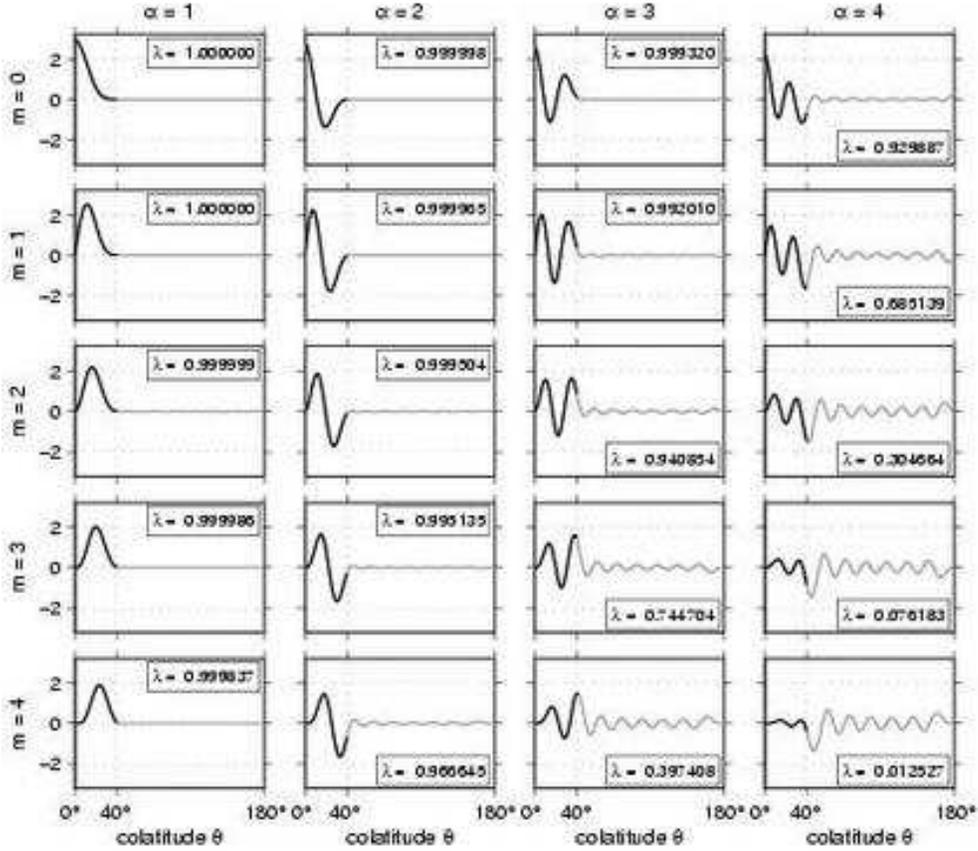}
} 
\caption{\small Colatitudinal dependence of the first four
spatial-domain eigenfunctions
$g_1(\theta),g_2(\theta),g_3(\theta),g_4(\theta)$ of fixed order $m=0$
(top) to $m=4$ (bottom).  The radius of the polar cap is
$\Theta=40^{\circ}$, and the maximal spherical harmonic degree is
$L=18$.  Black curves show the concentration within the cap
$0\leq\theta\leq 40^{\circ}$; grey curves highlight the leakage into the
rest of the sphere, $40^{\circ}<\theta\leq 180^{\circ}$.  The
eigenvalues $\lambda_1,\lambda_2,\lambda_3,\lambda_4$ expressing the
goodness of the spatial concentration are indicated. The corresponding
leakage in the spectral domain is illustrated in
Figure~\ref{sdwspectral}.}
\label{sdwspace} 
\end{figure}

We shall henceforth drop the identifying subscript
and rewrite the fixed-order, spectral-domain algebraic eigenvalue
problem $\sfD_{\pm m}\sfg_{\pm m}=\lambda_{\pm m}\sfg_{\pm m}$
as, simply,
\begin{equation}
\sfD\sfg=\lambda\hspace{0.1em}\sfg.
\label{fixedmeqn}
\end{equation}
The $(L-m+1)\times (L-m+1)$ matrix $\sfD$ and the $(L-m+1)$-dimensional
column vector $\sfg\in{\mathcal S}_L$ are of the form
\begin{equation}
\sfD = \left(\begin{array}{ccc}
D_{mm} & \cdots & D_{mL}\\
\vdots && \vdots\\
D_{Lm} & \cdots & D_{LL}\end{array}\right)
\quad\mbox{and}\quad
\sfg = \left(\begin{array}{c}
g_{m} \\ \vdots \\ g_{L}
\end{array}\right),
\label{littleD&g}
\end{equation}
where, for a particular order $m$,
\begin{equation}
D_{ll'}=2\pi\int_{0}^{\Theta}X_{lm}X_{l'm}\sin\theta\,d\theta.
\label{kernel4}
\end{equation}
The integral in equation~(\ref{kernel4}) can be evaluated
with the aid of equations~(\ref{byerly})--(\ref{legprod}):
\ber
\lefteqn{\hspace*{-18.0ex}
D_{ll'}=(-1)^m\frac{\sqrt{(2l+1)(2l'+1)}}{2}\!\!\sum_{n=|l-l'|}^{l+l'}
\left(\!\begin{array}{ccc}
l & n & l' \\ 0 & 0 & 0\end{array}\!\right)
\left(\!\begin{array}{ccc}
l & n & l' \\ m & 0 & \!\!-m\end{array}\!\right)}
\nonumber \\
&&\times\big[P_{n-1}(\cos\Theta)-P_{n+1}(\cos\Theta)\big].
\label{kernel3j}
\eer

We rank order the distinct $L-m+1$ eigenvalues
$\lambda_1,\lambda_2,\ldots,\lambda_{L-m+1}$ obtained by solving the
fixed-order eigenvalue problem~(\ref{fixedmeqn}) so that
\begin{equation}
1>\lambda_1>\lambda_2>\cdots>\lambda_{L-m+1}>0,
\label{fixedmorder}
\end{equation}
and orthonormalize the $(L-m+1)$-dimensional eigenvectors
$\sfg_1,\sfg_2,\ldots,\sfg_{L-m+1}$ as
\begin{equation}
\sfg_{\alpha}^{\sf{\sst{T}}}\sfg_{\beta}=\delta_{\alpha\beta}
\quad\mbox{and}\quad
\sfg_{\alpha}^{\sf{\sst{T}}}\sfD\hspace{0.05em}\sfg_{\beta}
=\lambda_{\alpha}\delta_{\alpha\beta}.
\label{1Dnormalize}
\end{equation}
The associated bandlimited eigenfunctions
$g_1(\theta),g_2(\theta),\ldots,g_{L-m+1}(\theta)$,
defined by
\begin{equation}
g(\theta)=\sum_{l=m}^Lg_lX_{lm}(\theta),
\label{polarg}
\end{equation}
then satisfy the colatitudinal orthogonality relations
\begin{equation}
2\pi\int_0^{\pi}g_{\alpha}g_{\beta}\sin\theta\,d\theta=\delta_{\alpha\beta}
\quad\mbox{and}\quad 2\pi\int_0^{\Theta}g_{\alpha}g_{\beta}\sin\theta\,d\theta
=\lambda_{\alpha}\delta_{\alpha\beta}.
\label{fixedmortho}
\end{equation}
The optimally concentrated spatial eigenfunctions $g(\rhat)$
for a given $m$ are expressed in terms of the fixed-order
colatitudinal eigenfunctions~(\ref{polarg}) by
\begin{equation}
g(\theta,\phi)=\left\{
\begin{array}{l@{\quad\mbox{if}\hspace{0.6em}}l}
\rule[-2mm]{0mm}{6mm}\sqrt{2}\,g(\theta)\cos m\phi & -L\le m<0\\
\rule[-2mm]{0mm}{6mm}g(\theta)                     & m=0\\
\rule[-2mm]{0mm}{6mm}\sqrt{2}\,g(\theta)\sin m\phi & 0< m\le L.\\
\end{array}
\right.
\label{polarg2}
\end{equation}

The four most optimally concentrated eigenfunctions
$g_1(\theta),g_2(\theta),g_3(\theta),g_4(\theta)$, for orders $0\leq
m\leq 4$ are plotted in Figure~\ref{sdwspace}. The radius of the polar
cap in this example is $\Theta=40^{\circ}$, and the maximal spherical
harmonic degree is $L=18$.  The first zonal ($m=0$) eigenfunction,
$g_1(\theta)$, has no nodes within the cap $0^{\circ}\leq\Theta\leq
40^{\circ}$; the second, $g_2(\theta)$, has one node, and so on.  The
non-zonal ($m>0$) eigenfunctions all vanish at the north pole,
$\Theta=0^{\circ}$.  The first four zonal eigenfunctions, the first
three $m=1$ and $m=2$ eigenfunctions, and the first two $m=3$ and
$m=4$ eigenfunctions are all very well concentrated ($\lambda>0.9$),
whereas the fourth $m=3$ and $m=4$ eigenfunctions exhibit significant
leakage ($\lambda<0.1$). The numerical methods used to compute the
results shown here and in ensuing figures are summarized in
Appendix~\ref{sec:Computational}. 

\subsection{Decomposition of the spatial-domain eigenvalue problem}
The integral eigenvalue problem~(\ref{firsttimeint}) in the spatial domain
likewise decomposes into a series of fixed-order,
one-dimensional Fredholm eigenvalue equations,
\begin{equation}
\int_{0}^{\Theta}D(\theta,\theta')\,g(\theta')\sin\theta'\,d\theta'
={\lambda}\hspace{0.1em}g(\theta),\quad 0\le\theta\le\pi,
\label{Fredholm}
\end{equation}
each with an $m$-dependent, separable, symmetric kernel,
\begin{equation}
D(\theta,\theta')=2\pi\suml_{l=m}^{L}X_{lm}(\theta)X_{lm}(\theta'). 
\label{dlx}
\end{equation}
\begin{figure}[b]\centering 
\rotatebox{0}{
\includegraphics[width=1\textwidth]{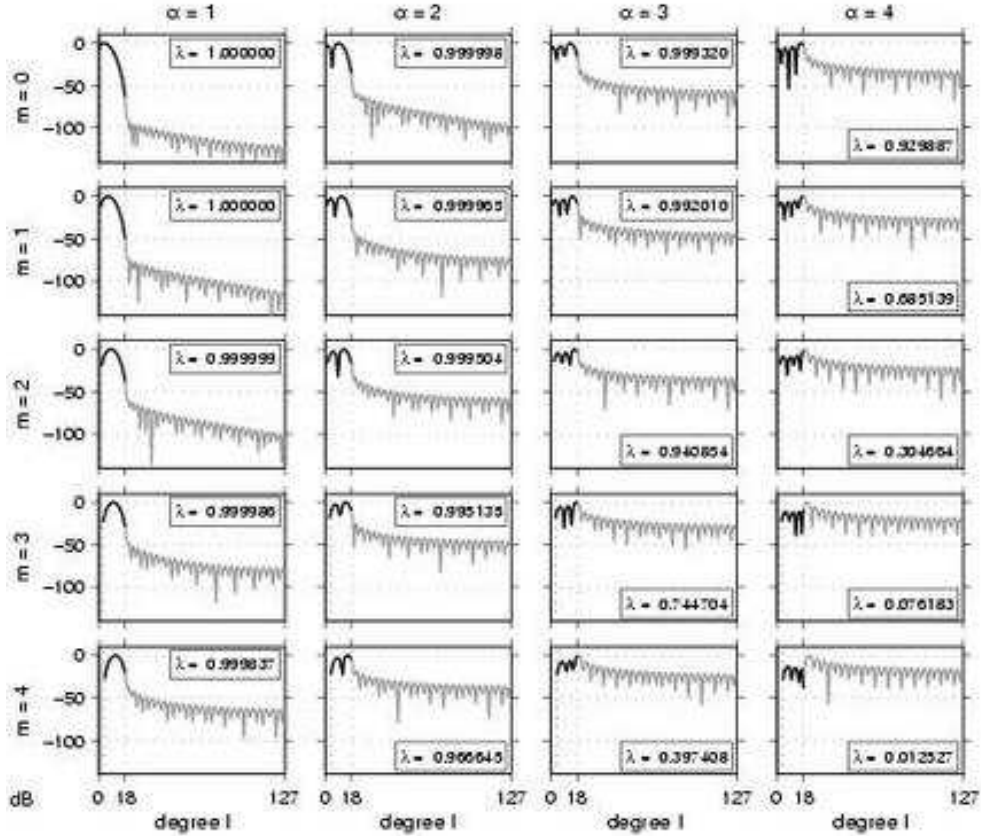}
} 
\caption{\small Squared coefficients $h_{l}^2$ of the first four
spacelimited eigenfunctions $h_1(\theta)$, $h_2(\theta)$,
$h_3(\theta)$, $h_4(\theta)$ of fixed order $m=0$
(top) to $m=4$ (bottom).  The radius of the polar cap is
$\Theta=40^{\circ}$, and the maximal spherical harmonic degree is
$L=18$.  Black curves show the power within the interval of
concentration $0\leq l\leq 18$; grey curves show the leaked power
within $19\leq l\leq 127$. Values of $h_{l}^2$ are in dB, normalized
to zero at the individual maxima. The eigenvalues
$\lambda_1,\lambda_2,\lambda_3,\lambda_4$ specifying the goodness of
the spectral concentration are indicated.  The corresponding
bandlimited, spatial-domain eigenfunctions are shown in
Figure~\ref{sdwspace}.}
\label{sdwspectral} 
\end{figure}
The results~(\ref{Fredholm}) and~(\ref{dlx}) can either be
obtained by multiplying the index form $D_{ll'}
g_{l'}=\lambda\hspace{0.1em}g_{l}$ of equation~(\ref{fixedmeqn})
by $X_{lm}(\theta)$ and summing over $m\le l\le L$,
or by substituting the representation~(\ref{polarg})--(\ref{polarg2})
into equation~(\ref{firsttimeint}), and using the orthogonality
of the longitudinal functions $\ldots,\sqrt{2}\cos m\phi,\ldots, 1,\ldots,
\sqrt{2}\sin m\phi,\ldots$ over the interval $0\le\phi< 2\pi$.
The $(L-m+1)\times (L-m+1)$ matrix eigenvalue problem~(\ref{fixedmeqn})
can in turn be derived starting from the separable Fredholm
equation~(\ref{Fredholm}), so that the fixed-order spectral and spatial
eigenvalue problems are completely equivalent. The fixed-order,
spacelimited eigenfunctions
\begin{equation}
h(\theta)=\left\{\begin{array}{ll}
g(\theta) & \mbox{if $0\le\theta\le\Theta$}\\
0 & \mbox{otherwise}
\end{array}\right.
\label{caphfromg}
\end{equation}
satisfy an equation identical to~(\ref{Fredholm}),
but only within the polar cap itself:
\begin{equation}
\int_{0}^{\Theta}D(\theta,\theta')\,h(\theta')\sin\theta'\,d\theta'
={\lambda}\hspace{0.1em}h(\theta),\quad 0\le\theta\le\Theta.
\label{Hredholm}
\end{equation}
The eigenvalue $\lambda$ is a measure
of both the spatial concentration of $g(\theta)\in{\mathcal S}_L$
within $0\le\theta\le\Theta$ and the spectral concentration
of $h(\theta)\in{\mathcal S}_R$ within $0\le l\le L$.
The substitution $\mu=\cos\theta$ converts
equations~(\ref{Fredholm}) and~(\ref{Hredholm}) into
\begin{subequations}
\label{bothholmmu}
\ber
\int_{\cos\Theta}^1D(\mu,\mu')\,g(\mu')\,d\mu'&=&
\lambda\hspace{0.1em}g(\mu), \quad-1\le\mu\le 1,\\
\label{Fredholmmu}
\int_{\cos\Theta}^1D(\mu,\mu')\,h(\mu')\,d\mu'&=&
\lambda\hspace{0.1em}h(\mu),\quad\cos\Theta\le\mu\le 1.
\label{Hredholmmu}
\eer
\end{subequations}
The kernel $D(\mu,\mu')$ can be simplified using
the Christoffel-Darboux identity~(\ref{christdarboux}):
\be
D(\mu,\mu')=\fracd{(L-m+1)!}{2(L+m)!}
\left[\fracd{P_{L+1,m}(\mu)P_{Lm}(\mu')-P_{Lm}(\mu)P_{L+1,m}(\mu')}
{\mu-\mu'}\right],
\label{CDkernel}
\ee where L'H\^opital's rule covers the case
$\mu=\mu'$. The squared spherical harmonic coefficients $h_{l}^2$ of
the four best concentrated spacelimited eigenfunctions
$h_1(\theta)$, $h_2(\theta)$, $h_3(\theta)$, $h_4(\theta)$ for $0\leq m\leq 4$
are plotted versus $l$ in Figure~\ref{sdwspectral}.  The cap radius
$\Theta=40^{\circ}$, bandwidth $L=18$, and layout are the same as in
Figure~\ref{sdwspace}. The maximal contribution to the $\alpha$th
zonal ($m=0$) eigenfunction comes from the harmonic degree satisfying
$\sqrt{l(l+1)}\approx \pi/(\alpha\Theta)$; physically, this
corresponds to fitting an integral number of asymptotic wavelengths
$2\pi/\sqrt{l(l+1)}$ within the cap of diameter $2\Theta$.
\begin{figure}[b]\centering 
\rotatebox{0}{
\includegraphics[width=0.85\textwidth]{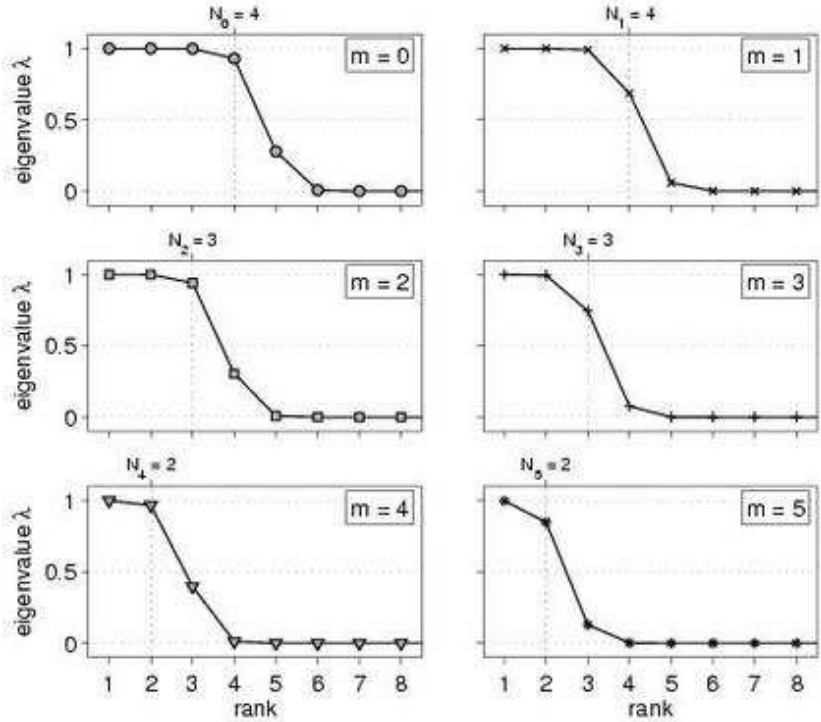}
} 
\caption{\small Fixed-order eigenvalue spectra for an axisymmetric
polar cap of radius $\Theta=40^{\circ}$.  The maximal spherical
harmonic degree is $L=18$.  A different symbol is used to plot
$\lambda_{\alpha}$ versus rank $\alpha$ for each order $0\leq m\leq
5$. The total number of fixed-order eigenvalues is $L-m+1$; only the
largest eight ($\lambda_1$ through $\lambda_8$) are shown.  Vertical
gridlines and top labels specify the partial Shannon numbers $N_m$,
rounded to the nearest integer.}
\label{sdweigen} 
\end{figure}
\subsection{Significant fixed-order eigenvalues}
The number of significant eigenvalues, or partial
Shannon number, for each of the fixed-order eigenvalue
problems~(\ref{fixedmeqn}), (\ref{Fredholm}),~(\ref{Hredholm})
or~(\ref{bothholmmu}) can be computed using any of the equivalent
formulae
\begin{equation}
N_m=\sum_{\alpha=1}^{L-m+1}\!\!\lambda_{\alpha}
=\suml_{l=m}^{L} D_{ll}
=\int_{0}^{\Theta}D(\theta,\theta)\,d\theta
=\int_{\cos\Theta}^1D(\mu,\mu)\,d\mu.
\label{nsubm}
\end{equation}
We can write the final relation in equation~(\ref{nsubm}) using
equation~(\ref{CDkernel}) in the form  
\begin{equation}
N_m=\fracd{(L-m+1)!}{2(L+m)!}\int_{\cos\Theta}^1
\big[P_{L+1,m}^{\prime}P_{Lm}-P_{Lm}^{\prime}P_{L+1,m}\big]\,d\mu,
\label{nsubm2}
\end{equation}
where the prime denotes differentiation with respect to $\mu$.
Figure~\ref{sdweigen} shows the fixed-order eigenvalue spectra for
$0\leq m\leq 5$. The cap radius is $\Theta=40^{\circ}$ and the maximal
spherical harmonic degree is $L=18$, as in Figures~\ref{sdwspace}
and~\ref{sdwspectral}. The partial Shannon numbers $N_m$, computed by
rounding equation~(\ref{nsubm}) to the nearest integer, are shown.  As
in the case of the classical Slepian problem
\cite[]{Landau65,Percival+93,Slepian+65} the spectra have a
characteristic step shape, showing significant ($\lambda\approx 1$)
and insignificant ($\lambda\approx 0$) eigenvalues separated by a
narrow transition band. The partial Shannon
number~(\ref{nsubm}) provides a good estimate of the number of
well concentrated eigenfunctions; the first $N_m$ eigenfunctions all
have a concentration factor exceeding $\lambda=0.5$.
\begin{figure}[b]\centering 
\rotatebox{0}{
\includegraphics[width=\textwidth]{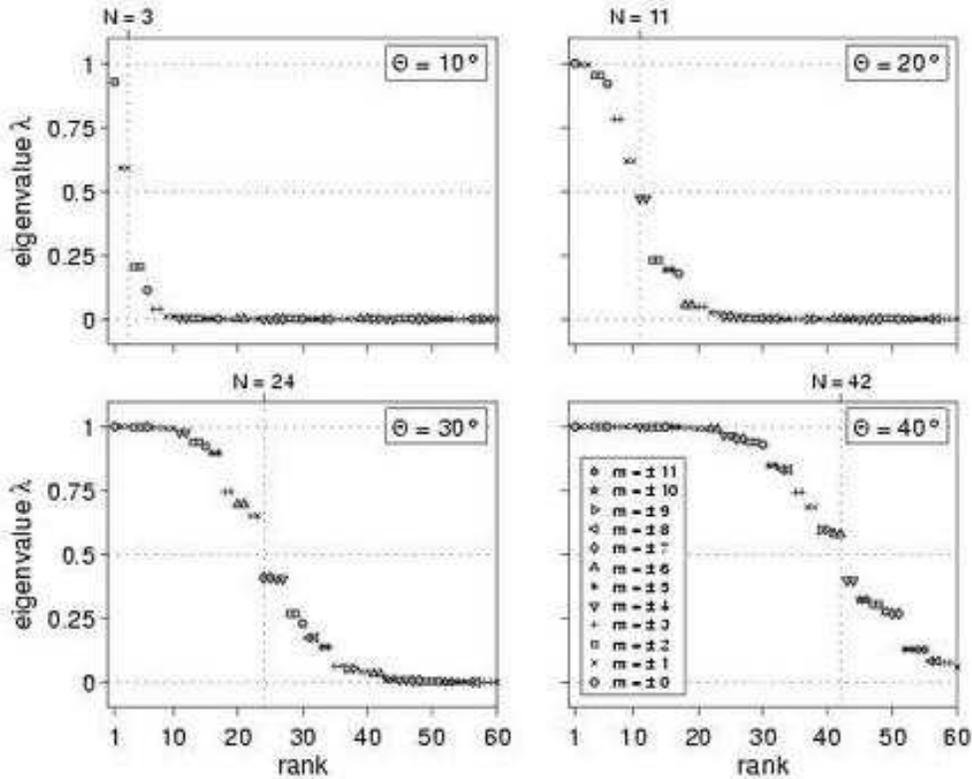}
} 
\caption{\small Reordered eigenvalue spectra ($\lambda_{\alpha}$
versus rank $\alpha$) for axisymmetric polar caps of
$\Theta=10^{\circ},20^{\circ},30^{\circ},40^{\circ}$ and
a common maximal spherical harmonic degree $L=18$. The total
number of eigenvalues is $\Lpot=361$;
only $\lambda_1$ through $\lambda_{60}$ are shown.
Different symbols are used to plot the various orders $-11\leq m\leq 11$;
the first six symbols are the same as those used in Figure~\ref{sdweigen}.
Vertical gridlines and top labels specify the rounded Shannon numbers.}
\label{sdwvals} 
\end{figure}

Once the $L+1$ sequences of fixed-order eigenvalues~(\ref{fixedmorder})
have been found, they can be resorted in accordance with the mixed-order
ranking~(\ref{eigorder}). The total number of significant eigenvalues is
\begin{equation}
\label{Nequalsum}
N=N_0+2\suml_{m=1}^{L}N_m,
\end{equation}
where the factor of two accounts for the $\pm m$ degeneracy. In
Figure~\ref{sdwvals} we show the reordered, mixed-$m$ eigenvalue
spectra for four different polar caps, with colatitudinal radii
$\Theta=10^{\circ},20^{\circ},30^{\circ},40^{\circ}$; the maximal
spherical harmonic degree is $L=18$. The rounded Shannon numbers
$N=3,11,24,42$, lie in the middle of the steep, transitional part of
the spectra, roughly separating the reasonably well concentrated
eigensolutions ($\lambda>0.5$) from the more poorly concentrated ones
($\lambda<0.5$) in all four cases.
\begin{figure}[b]\centering 
\rotatebox{0}{
\includegraphics[width=\textwidth]{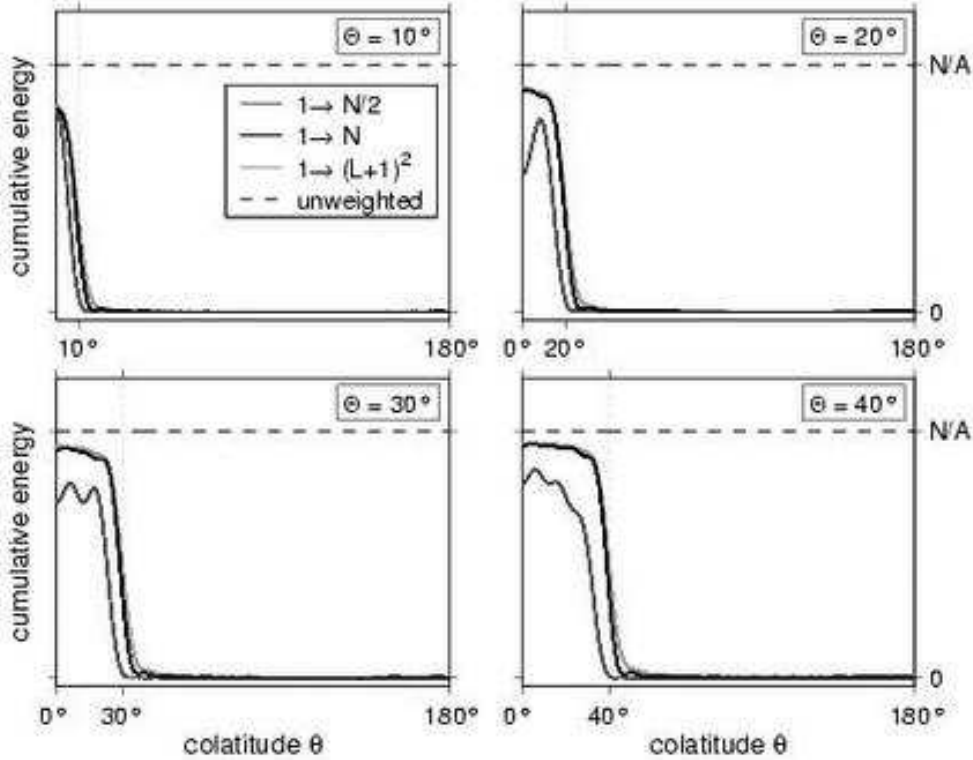}
} 
\caption{\small Cumulative energy of the eigenfunctions
concentrated within circularly symmetric polar caps of colatitudinal radii
$\Theta=10^{\circ},20^{\circ},30^{\circ},40^{\circ}$.
The maximal spherical harmonic degree is $L=18$;
the Shannon numbers are $N=3,11,24,42$. The sums
of squares $g_1^2(\theta,\phi)+g_1^2(\theta,\phi)+\cdots$
and $\lambda_1g_1^2(\theta,\phi)+\lambda_2g_2^2(\theta,\phi)+\cdots$,
are plotted versus colatitude $\theta$ along a fixed
arbitrary meridian $\phi$. Dashed lines show the full unweighted
sums of $\Lpot$ terms, which converge to the constant
value $N/A$ over the entire sphere $0\leq\theta\leq\pi$. Solid lines
show the eigenvalue-weighted partial sums of $N/2$ and $N$ terms,
and the full sum of $\Lpot$ terms, which are concentrated
within the polar cap $0\leq\theta\leq\Theta$.} 
\label{sdwsumall} 
\end{figure}

\begin{figure}[h]\centering 
\rotatebox{0}{
\includegraphics[width=1\textwidth]{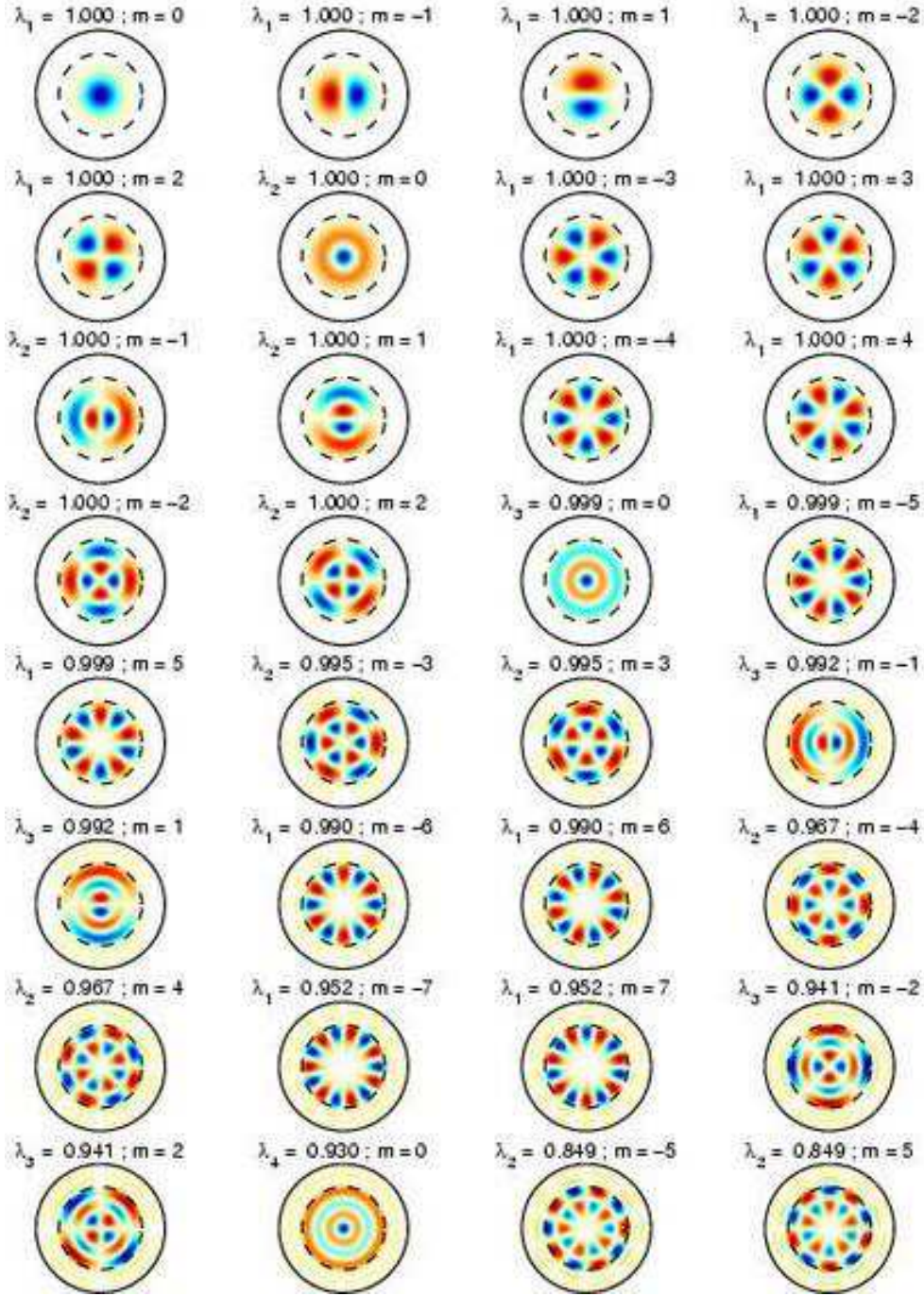}
} 
\caption{\small Bandlimited eigenfunctions $g(\theta,\phi)$
that are optimally concentrated within a circular cap of
 radius $\Theta=40^{\circ}$. Dashed circle
denotes the cap boundary. The maximal spherical harmonic
degree is $L=18$, and the Shannon number is $N=42$.  
Subscripts on the eigenvalues $\lambda_{\alpha}$ specify the
fixed-order rank. The eigenvalues have been resorted into
a mixed-order ranking, with the best concentrated eigenfunction
plotted on the top left and the 32nd best on the lower right.
Regions in which the absolute value is less than one
hundredth of the maximal absolute value on the sphere are
rendered in white; blue is positive, and red is negative.} 
\label{sdwfried} 
\end{figure}

We illustrate the pointwise sums of squares
$\sum_{\alpha}g_{\alpha}^2(\theta,\phi)$ and
$\sum_{\alpha}\lambda_{\alpha}g_{\alpha}^2(\theta,\phi)$ in
Figure~\ref{sdwsumall}, for polar caps of radii
$\Theta=10^{\circ},20^{\circ},30^{\circ},40^{\circ}$ and a bandwidth
$L=18$. The full unweighted sum of $\Lpot$ terms (dashed) is equal to
$N/A=\Lpot/(4\pi)$ over the entire sphere $0^{\circ}\leq\theta\leq
180^{\circ}$, in accordance with equation~(\ref{sumofsq}).  The
eigenvalue-weighted sums are, in contrast, concentrated within the
polar cap $0\leq\theta\leq\Theta$; solid \clearpage\noindent lines of different shades of
grey distinguish the sums carried up to the first $N/2$, $N$, or all
$(L+1)^2$ possible terms. Our expectation in equation~(\ref{sumofsq2})
is seen to be well realized. The uniformity of coverage over the
region of concentration achieved using only the first $N$ eigentapers
is a key result, responsible for the success of the multitaper method
of spectral estimation \cite[]{Walden90a}. The utilization of a
single, non-oscillatory data taper generally discards information near
the boundary; however, this information is recovered by application of
the other orthogonal tapers. The energy of the first $N$ tapered data
sets is proportional to the original energy of the data within the
region of concentration, enabling spatially uniform extraction of
statistical information while minimizing spectral leakage
\cite[]{Wieczorek+2004}. {\bf Mark's change here.}

Finally, Figure~\ref{sdwfried} shows a polar plot of the 32 eigenfunctions
$g(\theta,\phi)$, defined by equations~(\ref{polarg2}), that are optimally
concentrated within a polar cap of radius $\Theta=40^{\circ}$.  The
maximal spherical harmonic degree is $L=18$ and the Shannon number is
$N=42$, as in Figures~\ref{sdwspace}--\ref{sdweigen}.  The eigenvalue
ranking is mixed-order, as in Figure~\ref{sdwvals}, and all degenerate
$\sqrt{2}\cos m\phi, \sqrt{2}\sin m\phi$ doublets are shown. The
concentration factors $1<\lambda<0.849$ and orders $-5\leq m\leq 5$ of
each eigenfunction are indicated. Blue and red colors represent
positive and negative values, respectively; however, % the sign of every
%$g(\theta,\phi)$ is fundamentally arbitrary, and 
all signs could be
reversed without violating the quadratic concentration
criteria~(\ref{normratio}) and~(\ref{normratio2}).

\subsection{Commuting differential operator}
The analysis of the one-dimensional time-frequency
concentration problem was advanced considerably
by Slepian's seren\-dipitous discovery of the commuting prolate spheroidal
differential operator~(\ref{slepian7}). Remarkably, there is also
a differential operator, discovered by Gr\"{u}nbaum {\it et al.\/}
\cite[]{Grunbaum+82}, that commutes with the integral operator on the left side
of equations~(\ref{Hredholm}) and~(\ref{Hredholmmu}) \cite[see
also][]{Gilbert+77}. This second-order differential operator is
given explicitly by 
\begin{equation}
\label{grunbaumop}
\ssG=(\cos\Theta-\cos\theta)\nabla^2
+\sin\theta\,\frac{d}{d\theta}-L(L+2)\cos\theta,
\end{equation}
where
\begin{equation}
\nabla^2=\frac{d^2}{d\theta^2}+\cot\theta\,\frac{d}{d\theta}
-m^2(\sin\theta)^{-2}\label{fixedmlap1}
\end{equation}
is the fixed-order Laplace-Beltrami operator~(\ref{lapbelop}).
Rewritten in terms of $\mu=\cos\theta$,
the Gr\"{u}nbaum operator (\ref{grunbaumop}) is
\begin{equation}
\label{grunbaumop2}
\ssG=\fracd{d}{d\mu}\left[(\cos\Theta-\mu)
(1-\mu^2)\fracd{d}{d\mu}\right]-L(L+2)\mu
-\fracd{m^2(\cos\Theta-\mu)}{1-\mu^2}.
\end{equation}
Since commuting operators have the same eigenfunctions,
we can find the spacelimited, fixed-order eigenfunctions
$h(\theta)$ or $h(\mu)$ by solving the differential eigenvalue
equation $\ssG h=\chi\hspace{0.1em}h$, where $\chi\not= \lambda$ is the
associated Gr\"{u}nbaum eigenvalue.

To confirm that $\ssG$ in equation~(\ref{grunbaumop2})
is the desired commuting differential operator,
we are required to show that
\begin{equation}
\int_{\cos\Theta}^1\ssG_{\mu}D(\mu,\mu')\,h(\mu')\,d\mu'
=\int_{\cos\Theta}^1D(\mu,\mu')\,\ssG_{\mu'}\hspace{-0.01em}h(\mu')\,d\mu',
\label{grunshow0}
\end{equation}
for an arbitrary spacelimited function $h(\mu)$.
There are two steps in the argument given by Gr\"{u}nbaum {\it et al.\/}
\cite[]{Grunbaum+82}. The first is to show that the right
side of equation (\ref{grunshow0}) can be rewritten as
\begin{equation}
\int_{\cos\Theta}^1D(\mu,\mu')\,\ssG_{\mu'\hspace{-0.01em}}h(\mu')\,d\mu'
=\int_{\cos\Theta}^1\ssG_{\mu'}\hspace{-0.01em}D(\mu,\mu')\,h(\mu')\,d\mu',
\label{grunshow1}
\end{equation}
and the second is to show that
\begin{equation}
\ssG_{\mu}D(\mu,\mu')=\ssG_{\mu'}\hspace{-0.01em}D(\mu,\mu').
\label{grunshow2}
\end{equation}
The first result~(\ref{grunshow1}) is a straightforward exercise in integration
by parts; for any two functions $\zeta(\mu)$ and $\eta(\mu)$,
it may be easily shown that
\ber
\int_{\cos\Theta}^1\zeta(\ssG\eta)\,d\mu&=&
-\int_{\cos\Theta}^1\Big[(\cos\Theta-\mu)(1-\mu^2)\zeta'\eta'+L(L+2)\mu\,\zeta\eta
\nonumber \\
&&\qquad\qquad\,+\,m^2(\cos\Theta-\mu)(1-\mu^2)^{-1}\zeta\eta\Big]\,d\mu
=\int_{\cos\Theta}^1(\ssG\zeta)\eta\,d\mu,
\label{grunshow3}
\eer
where we have used a prime to denote $d/d\mu$.
Although it is not needed for the proof of equation~(\ref{grunshow0}),
we note for future reference that the result~(\ref{grunshow3}) is also
valid if the integrations are carried out over the full interval $-1\le\mu\le 1$.
To verify the second result~(\ref{grunshow2})
we make use of the relation $\nabla^2P_{lm}=-l(l+1)P_{lm}$
to write
\ber
\lefteqn{(\ssG_{\mu}-\ssG_{\mu'\hspace{-0.1em}})D(\mu,\mu')=}\hspace{3em}\nonumber \\
&&{}+\frac{1}{2}\sum_{l=m}^L\big[l(l+1)-L(L+2)\big]
(2l+1)\!\left[\fracd{(l-m)!}{(l+m)!}\right]\!P_{lm}(\mu)P_{lm}(\mu') \nonumber \\
&&{}-\frac{1}{2}(1-\mu^2)\sum_{l=m}^L(2l+1)\!\left[\fracd{(l-m)!}{(l+m)!}\right]\!
\frac{d}{d\mu}P_{lm}(\mu)P_{lm}(\mu') \nonumber \\
&&{}+\frac{1}{2}(1-\mu^{\prime 2})\sum_{l=m}^L(2l+1)
\!\left[\fracd{(l-m)!}{(l+m)!}\right]\!
P_{lm}(\mu)\frac{d}{d\mu'}P_{lm}(\mu').
\label{grunshow4}
\eer
An application of the Legendre derivative identity~(\ref{grunshow5})
transforms (\ref{grunshow4}) into
\ber
\lefteqn{
(\ssG_{\mu}-\ssG_{\mu'\hspace{-0.1em}})D(\mu,\mu')=}\hspace{3em}\nonumber \\
&&{}+\frac{1}{2}(\mu-\mu')\sum_{l=m}^L\big[l^2-(L+1)^2\big]
(2l+1)\!\left[\fracd{(l-m)!}{(l+m)!}\right]\!P_{lm}(\mu)P_{lm}(\mu') \nonumber \\
&&{}-\frac{1}{2}(\mu-\mu')\sum_{l=m}^L(2l+1)
\!\left[\fracd{(l-m+1)!}{(l+m)!}\right]\! \nonumber \\
&&{}\qquad\qquad\qquad\qquad\times
\big[P_{l+1\,m}(\mu)P_{lm}(\mu')-P_{lm}(\mu)
P_{l+1\,m}(\mu')\big],
\label{grunshow6}
\eer
and the Christoffel-Darboux identity~(\ref{christdarboux})
transforms equation~(\ref{grunshow6}) into
\ber
\lefteqn{
(\ssG_{\mu}-\ssG_{\mu'\hspace{-0.1em}})D(\mu,\mu')=}\hspace{3em}\nonumber \\
&&{}+\frac{1}{2}(\mu-\mu')\sum_{l=m}^L\big[l^2-(L+1)^2\big]
(2l+1)\!\left[\fracd{(l-m)!}{(l+m)!}\right]\!P_{lm}(\mu)P_{lm}(\mu') \nonumber \\
&&{}-\frac{1}{2}(\mu-\mu')\sum_{l=m}^L(2l+1)
\!\sum_{n=m}^l(2n+1)\!\left[\fracd{(n-m)!}{(n+m)!}\right]\!
P_{nm}(\mu)P_{nm}(\mu').
\label{grunshow7}
\eer
Interchanging the order of summation in the last line
of equation~(\ref{grunshow7}) we obtain\newpage
\ber
\lefteqn{
(\ssG_{\mu}-\ssG_{\mu'\hspace{-0.1em}})D(\mu,\mu')=}\hspace{3em} \nonumber \\
&&{}+\frac{1}{2}(\mu-\mu')\sum_{l=m}^L\big[l^2-(L+1)^2\big]
(2l+1)\!\left[\fracd{(l-m)!}{(l+m)!}\right]\!P_{lm}(\mu)P_{lm}(\mu') \nonumber \\
&&{}-\frac{1}{2}(\mu-\mu')\sum_{n=m}^L(2n+1)
\!\left[\fracd{(n-m)!}{(n+m)!}\right]\!
P_{nm}(\mu)P_{nm}(\mu')\sum_{l=n}^L(2l+1).
\label{grunshow8}
\eer
The final sum over $n$ is equal to $(L+1)^2-n^2$, so that the two
terms in equation~(\ref{grunshow8}) cancel,
$(\ssG_{\mu}-\ssG_{\mu'\hspace{-0.1em}})D(\mu,\mu')=0$, and the
commutation relation~(\ref{grunshow0}) is confirmed. 

\subsection{Gr\"{u}nbaum's equation}
The above argument shows that
we can compute the fixed-order, spacelimited,
colatitudinal eigenfunctions
$h_1(\theta),h_2(\theta),\ldots,h_{L-m+1}(\theta)$ either by solving
the integral equation~(\ref{Fredholm}) or by
solving the differential equation
\begin{equation}
(\cos\Theta-\cos\theta)\nabla^{2\hspace{-0.1em}} h
+\sin\theta\,\frac{dh}{d\theta}-L(L+2)\cos\theta\,h=\chi\hspace{0.1em}h,
\quad 0\le\theta\le\Theta.
\label{hrundiffeqn}
\end{equation}
The equivalent equation in terms of $\mu=\cos\theta$ is in
standard Sturm-Liouville form \cite[]{Courant+53}:
\begin{equation}
\label{sturm}
(ph')'-qh+\chi\rho h =0,
\quad \cos\Theta\le\mu\le 1,
\end{equation}
where the prime denotes differentiation with respect to $\mu$, and where
\begin{subequations}
\label{sturm2}
\ber
p(\mu)&=&(\cos\Theta-\mu)(1-\mu^2),\\
q(\mu)&=&m^2(1-\mu^2)^{-1}(\mu-\cos\Theta)-L(L+2)\mu,\\
\rho(\mu)&=&1.
\eer
\end{subequations}
Equation~(\ref{sturm}) must be solved subject to the requirement
that $h(\mu)$ remain finite at the endpoints $\mu=\cos\Theta$ and $\mu=1$.
The associated variational problem is \cite[]{Courant+53}
\begin{equation}
\label{sturm3}
\chi=\fracd{\int_{\cos\Theta}^1\big(ph^{\prime 2}+qh^2\big)\,d\mu}
{\int_{\cos\Theta}^1\rho h^2\,d\mu} =\mbox{minimum}.
\end{equation}
All of the familiar Sturm-Liouville theorems apply.  In particular, we
know that equation~(\ref{sturm}) has a simple spectrum, with an
infinite number of distinct eigenvalues $\chi_1<\chi_2<\ldots$, having
an accumulation point at infinity.  The rank orderings of the
eigenvalues $\chi_1,\chi_2,\ldots$ and the spatiospectral
concentration factors $\lambda_1,\lambda_2,\ldots,\lambda_{L-m+1}$ are
reversed, so that the eigenfunction $h_1(\theta)$ associated with the
numerically smallest eigenvalue $\chi_1$, which has no nodes in the
polar cap $0\le\theta\le\Theta$, is the most concentrated fixed-order
eigenfunction; $h_2(\theta)$, which has exactly one node, is the next
best concentrated, and so on.  Only the first $L-m+1$ eigenfunctions
$h_1(\theta),h_2(\theta),\ldots, h_{L-m+1}(\theta)$ with non-zero
eigenvalues $\lambda_1,\lambda_2,\ldots,\lambda_{L-m+1}$ are of interest
in most applications. The remaining eigenfunctions
$h_{L-m+2}(\theta),h_{L-m+3}(\theta),\ldots$ are in the null space of
the integral equation~(\ref{Hredholm}).

\subsection{Commuting tridiagonal matrix}
As in the case of~(\ref{Fredholm}) and~(\ref{Hredholm}),
we are free to extend the domain of equation~(\ref{hrundiffeqn})
to the entire domain $0\le\theta\le\pi$; in that case, the unknown function
must be bandlimited rather than spacelimited:
\begin{equation}
(\cos\Theta-\cos\theta)\nabla^{2\hspace{-0.1em}} g
+\sin\theta\,\frac{dg}{d\theta}-L(L+2)\cos\theta\,g=\chi\hspace{0.1em}g,
\quad 0\le\theta\le\pi.
\label{grundiffeqn}
\end{equation}
Upon substituting the harmonic
representation~(\ref{polarg}) of $g(\theta)$ into equation~(\ref{grundiffeqn}),
multiplying both sides by $2\pi\sin\theta\,X_{l'm}(\theta)$,
integrating over $0\le\theta\le\pi$, and invoking
the orthogonality relation~(\ref{Xlmortho}), we obtain
the algebraic eigenvalue equation
\begin{equation}
\sfG\hspace{0.05em}\sfg=\chi\hspace{0.1em}\sfg,
\label{grunmatrix}
\end{equation}
where
\begin{equation}
\sfG = \left(\begin{array}{ccc}
G_{mm} & \cdots & G_{mL}\\
\vdots && \vdots\\
G_{Lm} & \cdots & G_{LL}\end{array}\right)
\label{grunmatrix2}
\end{equation}
is the $(L-m+1)\times (L-m+1)$ matrix with elements
\begin{equation}
G_{ll'}=2\pi\int_{0}^{\pi}X_{lm}(\ssG\hspace{-0.05em}
X_{l'm})\sin\theta\,d\theta.
\label{grunmatrix3}
\end{equation}
Equation~(\ref{grunmatrix}) is the spectral-domain version
of the differential eigenvalue equation~(\ref{grundiffeqn})
just as equation~(\ref{fixedmeqn}) is the spectral-domain equivalent
of the integral eigenvalue equation~(\ref{Fredholm}).
\begin{figure}[t]\centering 
\rotatebox{0}{
\includegraphics[width=0.95\textwidth]{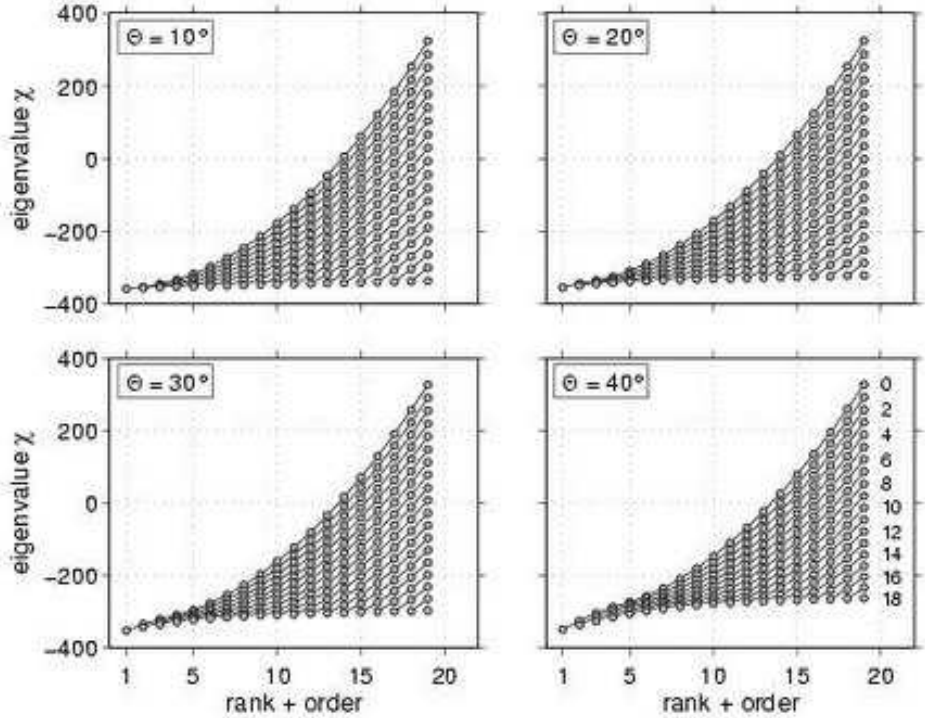}
} 
\caption{\small Rank-ordered Gr\"{u}nbaum eigenvalues for polar
caps of different colatitudinal radii,
$\Theta=10^{\circ},20^{\circ},30^{\circ},40^{\circ}$,
and a common maximal spherical harmonic degree $L=18$.
Separate sequences of eigenvalues $\chi_1,\chi_2,
\ldots,\chi_{L-m+1}$ for each angular order $0\le m\le L$ are connected
by lines, with each sequence offset horizontally by its order $m$ for clarity.}
\label{sdwgrunval} 
\end{figure}

As we have noted, the symmetry relation~(\ref{grunshow3})
is valid even if the interval of integration
is extended to $0\le\theta\le\pi$, as in
equation~(\ref{grunmatrix3}). This shows
that the Gr\"{u}nbaum matrix $\sfG$ is symmetric:
\begin{equation}
\sfG=\sfG^{\sf{\sst{T}}}.
\label{grunsymm}
\end{equation}
In addition, the matrices $\sfD$ and $\sfG$ commute with each other,
\begin{equation}
\sfD\hspace{0.05em}\sfG=\sfG\hspace{0.05em}\sfD,
\label{G&Dcommute}
\end{equation} 
so they have identical eigenvectors. The index version
of~(\ref{G&Dcommute}) is   
\begin{equation}
\sum_{n=m}^LD_{ln}G_{nl'}=2\pi\int_0^{\Theta}X_{lm}(\ssG
\hspace{-0.05em}X_{l'm})\sin\theta\,d\theta=\sum_{n=m}^LG_{ln}D_{nl'}.
\label{G&Dcommute2}
\end{equation}
The interior expression involving both the
operator $\ssG$ and integration over the region
of concentration $0\le\theta\le\Theta$ is the $ll'$ or $l'l$
element of the symmetric matrix product
$\sfD\hspace{0.05em}\sfG=(\sfD\hspace{0.05em}\sfG)^{\sf{\sst{T}}}$.
Verification of the intermediate steps requires
the use of the orthogonality relation~(\ref{Xlmortho}),
the operator identity~(\ref{grunshow2}),
and both the symmetry relation~(\ref{grunshow3}) and
its extension to the interval $0\le\theta\le\pi$.

There are a number of ways to evaluate the matrix
elements~(\ref{grunmatrix3}), perhaps the most straightforward of 
which is to make use of the relation $\nabla^2X_{lm}=-l(l+1)X_{lm}$
and the Legendre identities~(\ref{needrecur}). In fact, the
Gr\"{u}nbaum matrix $\sfG$ is tridiagonal: 
\begin{subequations}
\label{tridiag}
\ber
G_{ll}&=&-l(l+1)\cos\Theta,\label{tridiag1}\\
G_{l\,l+1} = G_{l+1\,l}&=&\big[l(l+2)-L(L+2)\big]
\sqrt{\fracd{(l+1)^2-m^2}{(2l+1)(2l+3)}},\label{tridiag2}\\
G_{ll'}&=&0\quad\mbox{otherwise},\label{tridiag3}
\eer
\end{subequations}
which is in agreement with the corresponding result given by
Gr\"{u}nbaum {\it et al.\/} \cite[]{Grunbaum+82}.
\newpage

The solution of equation~(\ref{grunmatrix}) offers a particularly
attractive means of computing the eigenvectors ${\textsf
g}\in{\mathcal S}_L$ and thus the optimally concentrated polar cap
eigenfunctions $g(\theta)\in{\mathcal S}_L$, because it only requires
the numerical diagonalization of a tridiagonal matrix $\sfG$ with
analytically prescribed elements~(\ref{tridiag}), and a spectrum of
eigenvalues $\chi$ that is guaranteed to be simple. To illustrate
this, we show the Gr\"{u}nbaum eigenvalue spectra for axisymmetric
polar caps of radii
$\Theta=10^{\circ},20^{\circ},30^{\circ},40^{\circ}$ and a maximal
spherical harmonic degree $L=18$ in Figure~\ref{sdwgrunval}.  The
ranked eigenvalues for every order $0\leq m\leq L$ are connected by
lines, each sequence offset by its order to facilitate
inspection.  Thus, $L+1$ eigenvalues $\chi_1,\chi_2,\ldots,\chi_{L+1}$
are plotted for $m=0$, whereas a single eigenvalue $\chi_1$ is plotted
for $m=L$. The roughly equant spacing and absence of a numerically
troublesome tail of near-zero values, as in
Figures~\ref{sdweigen}-\ref{sdwvals}, is evident.  

\begin{figure}[t]\centering 
\rotatebox{0}{
\includegraphics[width=0.9\textwidth]{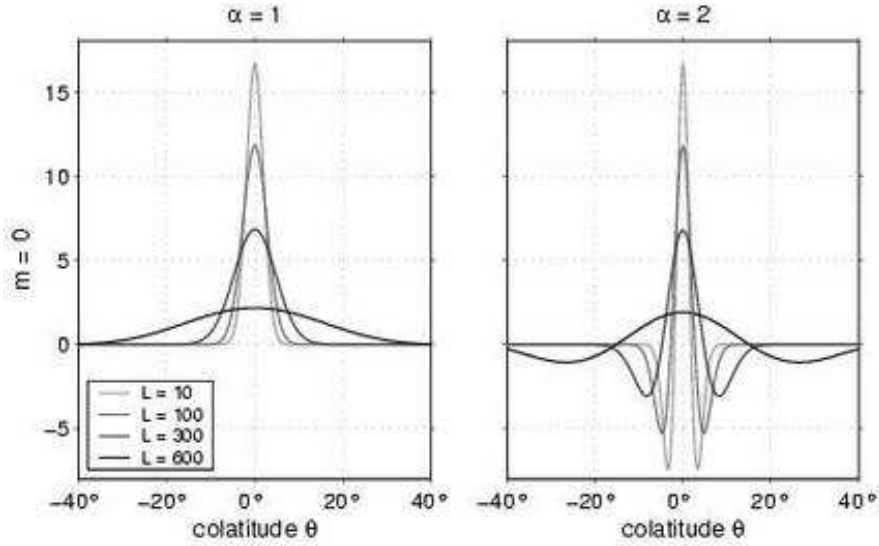}
} 
\caption{\small Zonal ($m=0$) eigenfunctions $g_1(\theta)$ (left) and
$g_2(\theta)$ (right) that are optimally concentrated within an
axisymmetric polar cap of colatitudinal radius $\Theta=40^{\circ}$,
for a range of maximal spherical harmonic degrees $L=10,100,300,600$
(increasing shades of grey). Roundoff error prevents the accurate
computation of the $L=300$ and $L=600$ eigenfunctions by
double-precision numerical diagonalization of the matrix $\sfD$;
diagonalization of the Gr\"{u}nbaum matrix $\sfG$ overcomes this
obstacle.}
\label{sdwspike} 
\end{figure}

Diagonalization of the Gr\"{u}nbaum matrix $\sfG$ enables the stable
computation of optimally concentrated spherical Slepian functions
that are bandlimited to very high angular degrees.
To illustrate this, we show in Figure~\ref{sdwspike} the first
two zonal ($m=0$) eigenfunctions $g_1(\theta),g_2(\theta)$
that are optimally concentrated within a polar cap of radius
$\Theta=40^{\circ}$, for increasing bandwidths $L=10,100,300,600$.
As the bandwidth increases, the eigenfunctions become increasingly
concentrated near the pole $\theta=0^{\circ}$. In any multitaper
spectral application, a greater and greater fraction of any data
near the boundary of the polar cap will be downweighted
upon windowing; however, with increasing bandwidth and thus Shannon number,
more and more well concentrated eigentapers become available,
to enable the recovery of this lost information.

\begin{figure}[h]\centering 
\rotatebox{0}{
\includegraphics[width=1\textwidth]{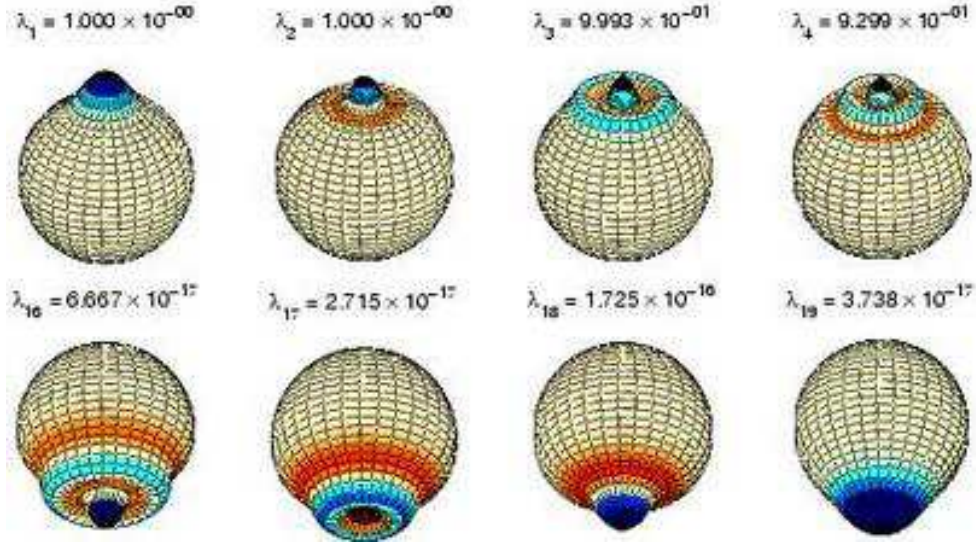}
} 
\caption{\small Optimally concentrated (top row) and optimally
excluded (bottom row) zonal ($m=0$) eigenfunctions,
for a circular polar cap of colatitudinal radius $\Theta=40^{\circ}$
and a maximal spherical harmonic degree $L=18$. The optimally
excluded eigenfunctions cannot be accurately computed
by double-precision diagonalization of the matrix $\sfD$.
Solution of the concentration problem for a polar cap
of radius $\Theta=140^{\circ}$ gives rise to the same eigenfunctions
$g_1(\theta),g_2(\theta),g_3(\theta),g_4(\theta),\ldots,
g_{16}(\theta),g_{17}(\theta),g_{18}(\theta),g_{19}(\theta)$,
but in reverse order.} 
\label{sdwcmb} 
\end{figure}

Concentration within a large rather than a small polar cap is another
problem for which Gr\"unbaum's method is ideally suited.  To
illustrate this, we show in Figure~\ref{sdwcmb} the first four
($g_1,g_2,g_3,g_4$) and the last four ($g_{16},g_{17},g_{18},g_{19}$)
zonal ($m=0$) eigenfunctions, once again for a polar cap of radius
$\Theta=40^{\circ}$ and a maximal spherical harmonic degree $L=18$.
As noted in Section~\ref{Section:sigeivs}, the eigenfunctions that are
optimally excluded from the polar cap $\Theta=40^{\circ}$ are
optimally concentrated within the much larger antipodal cap
$\Theta=140^{\circ}$. The actual eigenvalues
$\lambda_{16},\lambda_{17},\lambda_{18},\lambda_{19}$ are many orders
of magnitude smaller than the listed values, which simply represent
the noise floor of our double-precision computations. The optimally
excluded eigenfunctions can nevertheless be accurately computed by
diagonalizing the tridiagonal matrix $\sfG$; this is
not otherwise possible in double precision arithmetic (see
Appendix~\ref{sec:Computational}).

\subsection{Abstract operator formulation}
The spatial-domain commutation relation~(\ref{grunshow0})
can be expressed using the operator notation of
Section~\ref{subsec:AbsOps} as
\begin{equation}
(\ssR\ssH^{-1}\ssL\ssH\ssR)\ssG=\ssG(\ssR\ssH^{-1}\ssL\ssH\ssR).
\label{abgruncomm}
\end{equation}
Since the Gr\"{u}nbaum operator $\ssG$ acts only upon spacelimited
colatitudinal functions $h(\theta)$, it must satisfy
\begin{equation}
\ssG=\ssG\ssR=\ssR\ssG.
\label{abgrunspace}
\end{equation}
Upon pre-multiplying equation~(\ref{abgruncomm}) by $\ssL\ssH$,
post-multiplying it by $\ssH^{-1}\ssL$, and making use of
the relation~(\ref{abgrunspace}) and the fact that $\ssL^2=\ssL$,
we obtain
\begin{equation}
(\ssL\ssH\ssR\ssH^{-1}\ssL)(\ssL\ssH\ssG\ssH^{-1}\ssL)=
(\ssL\ssH\ssG\ssH^{-1}\ssL)(\ssL\ssH\ssR\ssH^{-1}\ssL).
\label{abgruncomm2}
\end{equation}
Equation~(\ref{abgruncomm2}) is the abstract operator
formulation of the spectral-domain matrix commutation relation
${\textsf D}\sfG=\sfG{\textsf D}$.
The operator equivalents of the spectral-domain
eigenvalue equation~(\ref{grunmatrix}) and the
spatial-domain eigenvalue equation~(\ref{hrundiffeqn}) are
\begin{subequations}
\label{finopeqns}
\ber
(\ssL\ssH\ssG\ssH^{-1}\ssL)(\ssL\,\sff)&=&
\chi\hspace{0.1em}(\ssL\,\sff),\\
\ssG(\ssR f)&=&\chi\hspace{0.1em}(\ssR f),
\eer
\end{subequations}
where $f(\theta)$ is an arbitrary colatitudinal
function, neither bandlimited nor spacelimited.
Because of the commutation relations~(\ref{abgruncomm})
and~(\ref{abgruncomm2}) we are free to solve
equations~(\ref{finopeqns}) rather than
the fixed-order version of equations~(\ref{gruncomp}),
to find the bandlimited eigenvectors
$\sfg=\ssL\,\sff$
and the spacelimited eigenfunctions
$h(\theta)=\ssR f(\theta)$.

\section{Continental Concentration}
\label{sec:Numerical}

\begin{figure}[b]\centering 
\rotatebox{0}{
\includegraphics[width=0.94\textwidth]{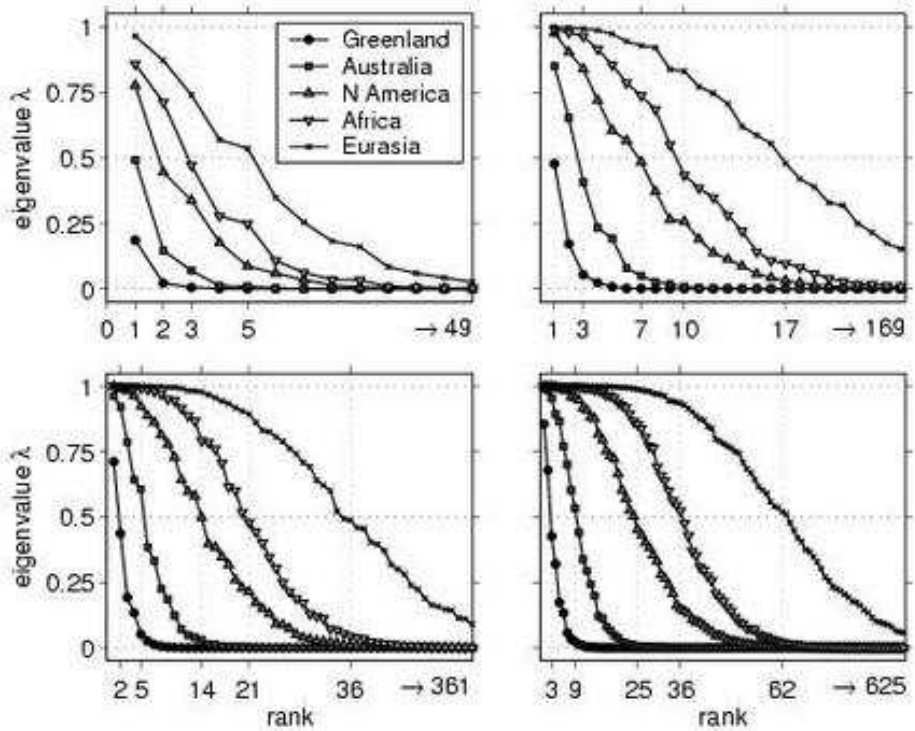}
} 
\caption{\small Eigenvalue spectra for five of the Earth's continental
regions (Greenland, Australia, North America, Africa, Eurasia) and
four different maximal spherical harmonic degrees ($L=6,12,18,24$).
Vertical gridlines and five leftmost ordinate labels specify the rounded Shannon
numbers $N$. Ordinates are truncated on the right; number to right of
arrow is the total number of eigenvalues, $\Lpot=49,169,361,625$.}
\label{sdwregions} 
\end{figure} 

To illustrate the theory for an irregularly shaped region,
we consider the spatiospectral concentration problem in
six of the Earth's continental regions, listed in Table~\ref{tebel}
in order of increasing size, together with their Shannon numbers for
different bandwidths. The spectral analysis of data within either
the Earth's continents or oceans has a number of
applications in geodesy, geophysics and oceanography
\cite[e.g.,][]{Hwang93,Hwang+97,Pail+2001,Shapiro+99b};
the spherical Slepian multitapers illustrated here should be ideally
suited for this task. 
\begin{table}[h]\centering{
\begin{tabular}{c|c|ccccc}\hline
Continental        & Area $A/(4\pi)$ & \multicolumn{4}{c}{Shannon
  number $N$}\\
region & in \% & $L=6$ & $L=12$ & $L=18$ & $L=24$  \\\hline
Greenland       & 0.44 & 0 & 1 &2 &3 \\
Australia       & 1.50 & 1 & 3 &5 &9 \\
South America   & 3.50 & 2 & 6 &13 &22 \\
North America   & 3.98 & 2 & 7 &14 &25 \\
Africa          & 5.78 & 3 & 10 &21 &36 \\
Eurasia         & 9.98 & 5 & 17 &36 &62\\\hline
\end{tabular}}\vspace{1em}\caption{\label{tebel}Areas, Shannon
numbers, and bandwidths for the continental concentration problem.}
\end{table}

Figure~\ref{sdwregions} shows the eigenvalue spectra for five of the
six regions (Greenland, Australia, North America, Africa and Eurasia)
and four different bandwidths, $L=6,12,18,24$,
corresponding to $\Lpot=49,169,361,625$ eigenfunctions each.  The
cutoff wavenumber associated with a bandwidth limit $L$ is
$\sqrt{L(L+1)}\approx L+1/2$ divided by the Earth's radius
\cite[]{Brune64,Jeans23}; the cutoff wavelengths corresponding to the
choices $L=6,12,18$ and $24$ are $6200,3200,2200$ and $1600$
kilometers, respectively. Only the largest of the Earth's continents,
Eurasia, is sizable enough to exhibit at least one nearly perfectly
concentrated eigenfunction for the smallest degree, $L=6$, and the
smallest region considered, Greenland, is too tiny to exhibit even a
single eigenfunction with a concentration factor $\lambda$ near unity
for the largest degree, $L=24$.  As in the case of a polar cap
(Figure~\ref{sdwvals}), the rounded Shannon numbers $N=\Lpot A/(4\pi)$
shown by the vertical dotted lines roughly separate the eigenfunctions
with concentration factors $\lambda>0.5$ from those with concentration
factors $\lambda>0.5$.

\begin{figure}[t]\centering 
\rotatebox{0}{
\includegraphics[width=1\textwidth]
{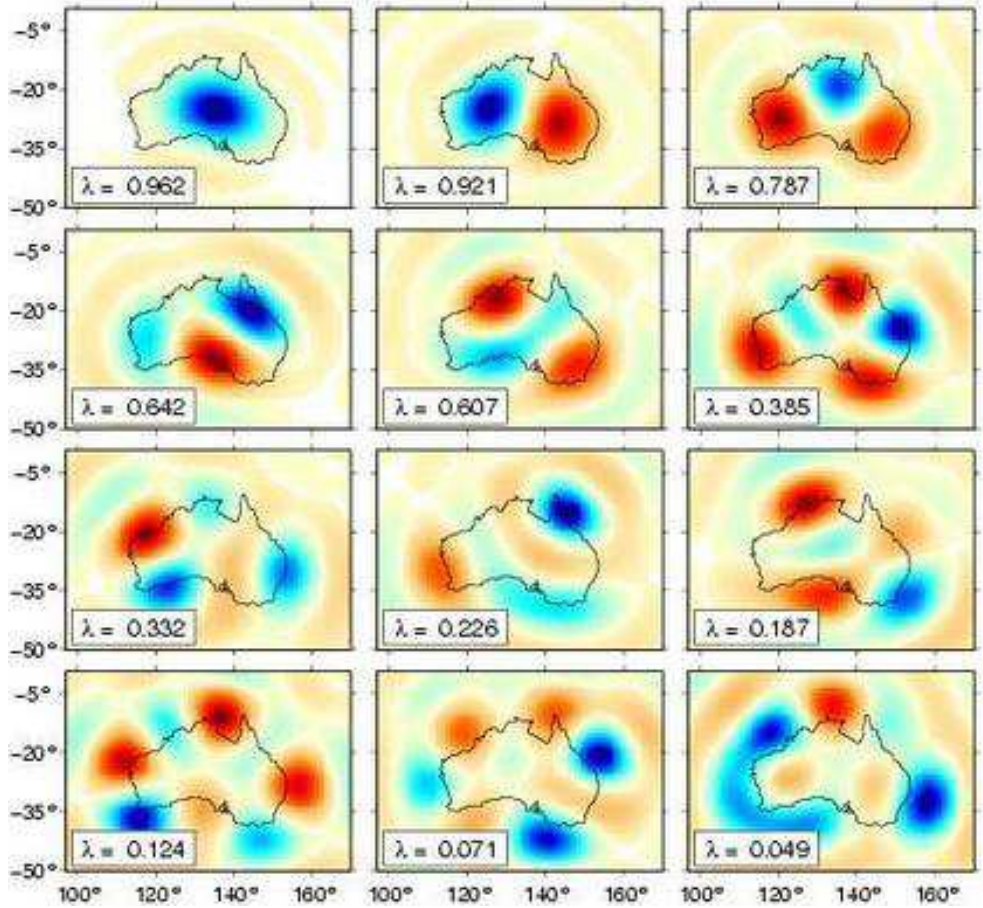}
} 
\caption{\small Bandlimited $L=18$ eigenfunctions $g_1,g_2,\ldots,g_{12}$
that are optimally concentrated within the continent of Australia. The
concentration factors $\lambda_1,\lambda_2,\ldots,\lambda_{12}$
are indicated; the Shannon number is $N=5$. Order is left to right,
top to bottom, as with English text.}
\label{sdwregions_australia}
\end{figure}

\begin{figure}[t]\centering 
\rotatebox{0}{
\includegraphics[width=1\textwidth]
{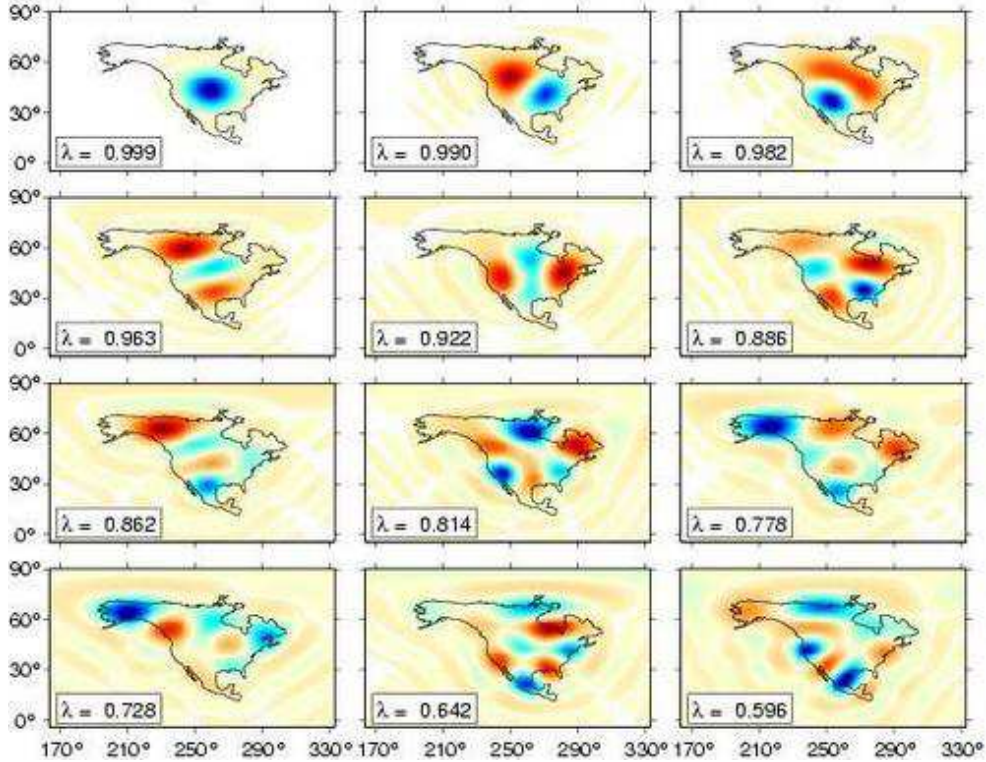}
} 
\caption{\small Bandlimited $L=18$ eigenfunctions $g_1,g_2,\ldots,g_{12}$
that are optimally concentrated within the continent of North America. The
concentration factors $\lambda_1,\lambda_2,\ldots,\lambda_{12}$
are indicated; the Shannon number is $N=14$. Format is identical to
that in Figure~\ref{sdwregions_australia}.}
\label{sdwregions_namerica} 
\end{figure}

\begin{figure}[t]\centering 
\rotatebox{0}{
\includegraphics[width=0.75\textwidth]
{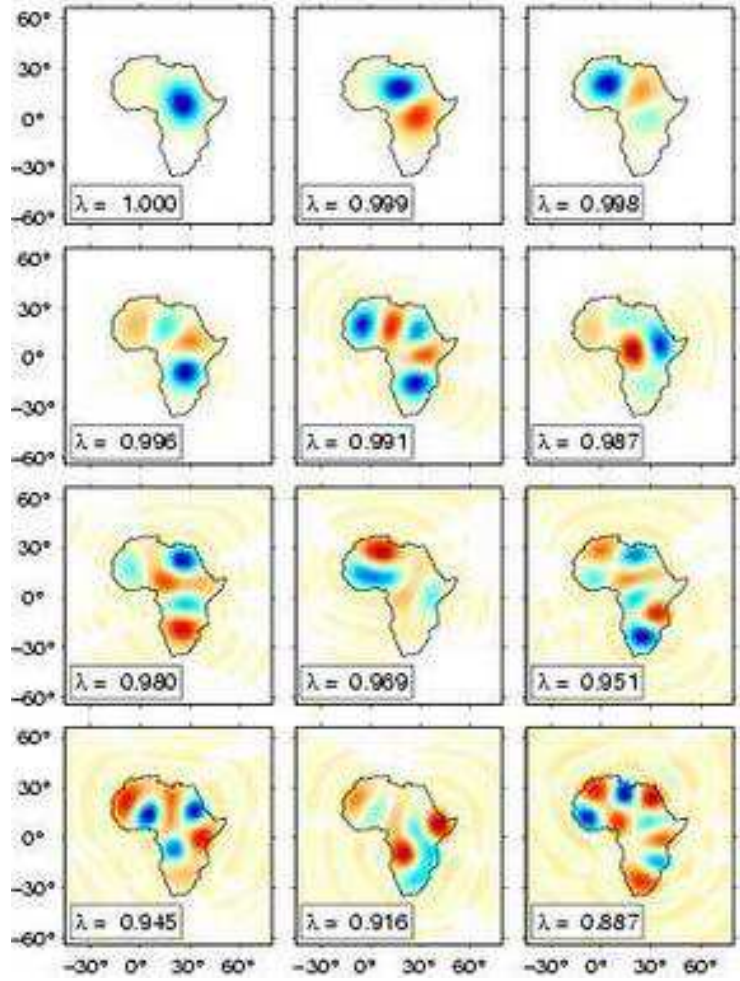}
} 
\caption{\small Bandlimited $L=18$ eigenfunctions $g_1,g_2,\ldots,g_{12}$
that are optimally concentrated within the continent of Africa. The
concentration factors $\lambda_1,\lambda_2,\ldots,\lambda_{12}$
are indicated; the Shannon number is $N=21$.  Format is identical to
that in Figure~\ref{sdwregions_australia}.}
\label{sdwregions_africa} 
\end{figure}

In Figures~\ref{sdwregions_australia},~\ref{sdwregions_namerica}
and~\ref{sdwregions_africa} we show map views of the first twelve
$L=18$ eigenfunctions $g_1(\rhat),g_2(\rhat),\ldots, g_{12}(\rhat)$
that are optimally concentrated within Australia, North America and
Africa.  Blue colors denote positive values and red colors denote
negative values (though, as we have noted, these could be reversed,
since the sign of an eigenfunction is arbitrary). Regions in which the
absolute value is less than one hundredth of the maximal absolute
value on the sphere are rendered in white. 

In the case of Australia (Figure~\ref{sdwregions_australia})
the first five eigenfunctions are reasonably well concentrated within
the continental boundaries ($\lambda_5=0.607$);
however, the concentration factors $\lambda$ diminish rapidly
thereafter, so that $g_{12}$ is far more excluded than concentrated
($\lambda_{12}=0.049$).  With a limiting bandwidth $L=18$,
and thus a cutoff wavelength of 2200 kilometers, it is only
possible to concentrate $N=5$ orthogonal bandlimited
eigenfunctions $g_1,g_2,g_3,g_4,g_5$ into a continent which,
across its north-south waist, is only about $1500$ kilometers wide.

This situation is much improved in the case of the North American
continent (Figure~\ref{sdwregions_namerica}), which has an area $A$
that is $2.7$ times larger than the area of Australia. In fact, North
America has $N=14$ reasonably well concentrated $L=18$ eigenfunctions,
of which only the first twelve are shown.  The first eigenfunction,
$g_1$, is a roughly circular dome centered in the middle of the
continent, as in the case of Australia.  Subsequent orthogonal
eigenfunctions $g_2,g_3,\ldots$ exhibit lobes in previously uncovered
regions.  The coastline of North America is more irregular than that
of Australia; Qu\'ebec and the Northwest Territories of Canada are
essentially uncovered until $g_8$, western Alaska is not adequately
covered until $g_9$ and $g_{10}$, and Florida, Baja California and
southern Mexico are only covered by $g_{11}$ and $g_{12}$ at the
expense of substantial leakage
($\lambda_{11}=0.642,\lambda_{12}=0.596$) outside of the continental
boundaries.

Africa (Figure~\ref{sdwregions_africa}), which has an area $A$ that is
3.9 times larger than that of Australia, has $N=21$ reasonably well
concentrated $L=18$ eigenfunctions, the twelfth of which has a
concentration factor $\lambda_{12}=0.887$.  Once again, $g_1$ is
roughly circular, and the subsequent orthogonal eigenfunctions
$g_2,g_3,\ldots$ successively cover previously uncovered regions.
West Africa is uncovered by $g_1$ and $g_2$, but becomes reasonably
well covered by $g_3$ and $g_5$; likewise, South Africa is uncovered
until $g_4$ and $g_5$. Other geographical features become well covered
by the increasingly oscillatory orthogonal eigenfunctions (e.g., Egypt
by $g_7$ and $g_{12}$).

\begin{figure}[h]\centering 
\rotatebox{0}{
\includegraphics[width=1\textwidth]{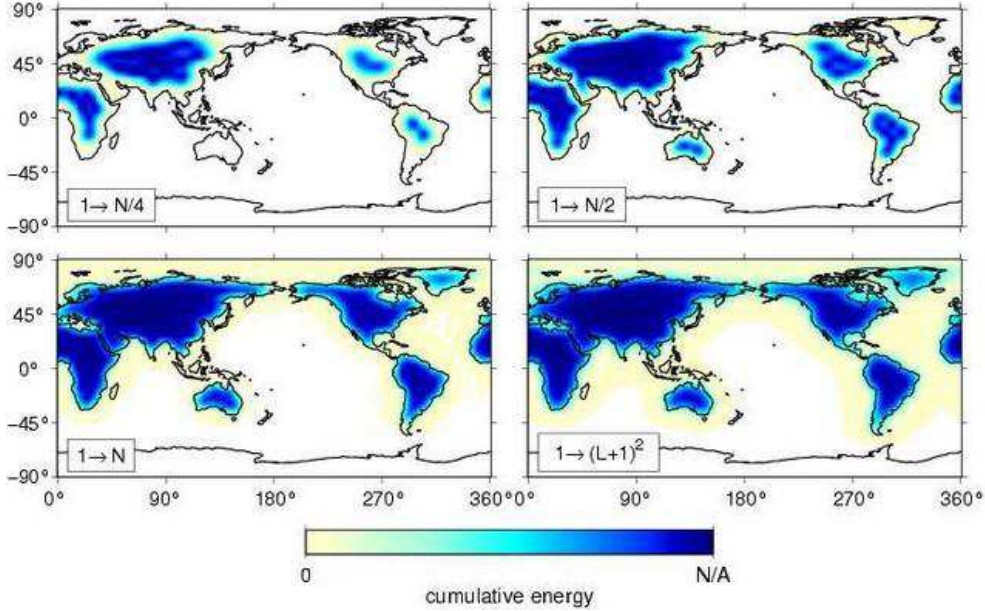}
} 
\caption{\small Cumulative eigenvalue-weighted energy
of the first $N/4,N/2,N$ and all $\Lpot$ eigenfunctions
that are optimally concentrated within the ensemble
of continents Eurasia, Africa, North America,
South America, Australia and Greenland. The maximal
spherical harmonic degree is $L=18$; the cumulative
fractional area is $A/(4\pi)=25.2\%$; the Shannon
number is $N=91$. Darkest blue on color bar corresponds
to the expected value~(\ref{sumofsq2}) of the sum,
as shown. Regions in which the value is less than one
hundredth of the maximal value on the sphere are
rendered in white.}
\label{sdwconsum} 
\end{figure}

Figure~\ref{sdwconsum} shows the eigenvalue-weighted sum of squares
$\sum_{\alpha}\lambda_{\alpha}g_{\alpha}^2(\rhat)$ of the bandlimited
$L=18$ eigenfunctions of all six of the Earth's continents (excluding
Antarctica).  We find the eigenfunctions
$g_1(\rhat),g_2(\rhat),\ldots,g_{\Lpot}(\rhat)$ by diagonalization of
the $\Lpot\times\Lpot$ matrix~(\ref{Dmatrixdef}) formed by summing the
corresponding matrices $\sfD_{\rm Eurasia}+\sfD_{\rm Africa}+ \cdots$ of
each of the six continents. The combined area of all six continents is
$A/(4\pi)=25.2\%$, and the Shannon number is $N=91$; the partial sums
of the first $N/4,N/2$ and $N$ terms, as well as the full sum of all
$\Lpot=361$ terms, are shown.  The ability of the first $N$
eigenfunctions to provide uniform coverage of the target area is
evident; as in Figure~\ref{sdwsumall}, the coverage is only marginally
improved by adding the remaining, poorly concentrated $\Lpot-N=250$
terms.  Due to their small size, Australia and Greenland do not appear
until the $1\rightarrow N/2$ and $1\rightarrow N$ partial sums,
respectively. Even then, the coverage of Greenland is imperfect, an
expected consequence of the small Shannon number for Greenland ($N=2$
for a maximal spherical harmonic degree $L=18$).

\newpage
\section{Asymptotic scaling}
\label{sec:Asymptotics}
As we have noted, the eigenvalues $\lambda_1,\lambda_2,\ldots$
and suitably scaled eigenfunctions $\psi_1(x),\psi(x)_2,\ldots$
of the original Slepian concentration problem~(\ref{slepian5})
depend only upon the Shannon number $N=2TW/\pi$. This
Shannon-number scaling is the only important feature of the
one-dimensional problem that does not carry over to the
spatiospectral concentration problem on a sphere.
Fundamentally, this lack of scaling is a consequence
of the fact that it is not possible to shrink or
magnify a region, such as Africa, on a sphere $\Omega$
of fixed radius $\|\rhat\|=1$, while keeping the angular
relationships among all of the interior points the same.
Shannon-number scaling on a sphere is exhibited only
asymptotically, in the limit
\begin{equation}
A\rightarrow 0,\quad L\rightarrow\infty,\quad
\mbox{with}\quad N=\Lpot\fracd{A}{4\pi} \quad\mbox{held fixed}.
\label{asylimit}
\end{equation}
In that limit of a small concentration area $A$ and a
large bandwidth $0\le l\le L$, the curvature of the
sphere becomes negligible and the spherical
concentration problem approaches the concentration problem
on a plane.  

\subsection{Hilb approximation and Poisson sum formula}
Two  results underlie the consideration of the flat-earth
limit~(\ref{asylimit}), which we undertake in this section.
The first is Hilb's asymptotic
approximation for the Legendre functions
\cite[]{Amado+85,Dahlen80,Hilb19,Szego75}, 
\begin{equation}
X_{lm}(\theta)\approx (-1)^m\sqrt{\frac{l+1/2}{2\pi}}\,
\sqrt{\frac{\theta}{\sin\theta}}\,
J_m\big[(l+1/2)\theta\big],\quad 0\le\theta\ll\pi,
\label{Hilb}
\end{equation}
where $J_m(x)$ is the Bessel function of the first kind,
and the second is the truncated Poisson sum formula,
\begin{equation}
\suml_{l=0}^L{f}(l+1/2)=\sum_{s=-\infty}^{\infty}(-1)^s
\int_{0}^{L+1}{f}(k)e^{-2\pi isk}\,dk,
\label{Poisson}
\end{equation}
which is valid for an arbitrary continuous function $f(x)$.
To verify the relation~(\ref{Poisson}) we start with
the Fourier series representation of $f(x)$
on the interval $0\le x\le 2\pi$,
\begin{equation}
f(x)=\frac{1}{2\pi}\sum_{s=-\infty}^{\infty}
\int_0^{2\pi}f(u)e^{is(x-u)}\,du.
\label{fishyproof1}
\end{equation}
Letting $x\rightarrow x+2\pi l$ in equation~(\ref{fishyproof1})
and summing over $0\le l\le L$ yields
\begin{equation}
\sum_{l=0}^Lf(x+2\pi l)=\frac{1}{2\pi}\sum_{s=-\infty}^{\infty}
\int_0^{(L+1)2\pi}f(u)e^{is(x-u)}\,du,
\label{fishyproof2}
\end{equation}
where we have shifted the interval of integration for each
term by $2\pi l$ and made use of the $2\pi$ periodicity
of the exponential. Upon dividing the argument by $2\pi$,
substituting $k=u/(2\pi)$, and setting $x=\pi$, we obtain
the desired identity~(\ref{Poisson}).

\subsection{Scaled integral equation for an arbitrary region}
An application of both the Hilb approximation~(\ref{Hilb}) and the
Poisson sum formula~(\ref{Poisson}) enables us to write the
Fredholm kernel $\Drhrhp$ in equation~(\ref{firsttimeint}) in the form
\ber
\DDlta&=&
\suml_{l=0}^{L}\tlofp\!P_l(\cos\Delta)\nnr\\
&\approx&\frac{1}{2\pi}\,
\sdsind \suml_{l=0}^{L}(l+1/2)J_0\big[(l+1/2)\Delta\big]\nnr\\
&=&\frac{1}{2\pi}\,\sdsind 
\suml_{s=-\infty}^{\infty}(-1)^s\int_{0}^{L+1}
J_0(k\Delta)e^{-2\pi isk}k\,dk.\label{kernel40}
\eer
Upon substituting $k=(L+1)p$ and taking
the limit $L\rightarrow\infty,\Delta\rightarrow 0$,
with the product $L\Delta$ held fixed, equation~(\ref{kernel40})
reduces to
\begin{equation}
D(\Delta)\approx\frac{\Lpot}{2\pi}
\int_{0}^{1}J_0\big[(L+1)p\Delta\big]\,p\,dp\nnr
=\frac{(L+1)\,J_1\big[(L+1)\Delta\big]}{2\pi\Delta},
\label{smallK}
\end{equation}
where we have made the approximation $\Delta/\!\sin\Delta\approx 1$,
and used the Riemann-Lebesgue lemma \cite[]{Olver97} to eliminate the
$s\not= 0$ terms involving the highly oscillatory factors $e^{-2\pi is(L+1)p}$.
In the limit $x\rightarrow 0$ the ratio $J_1(x)/x\rightarrow 1/2$,
so the $\Delta\rightarrow 0$ limit of the kernel~(\ref{smallK}) is $D(0)=\Lpot/(4\pi)$,
guaranteeing that the Shannon number, or sum of the eigenvalues~(\ref{tracedef}),
is still given in this asymptotic approximation by
\begin{equation}
\label{asyshannon}
N=\intr D(0)\domg=\Lpot\,\frac{A}{4\pi}.
\end{equation}

To obtain a scaled version of equation~(\ref{firsttimeint})
dependent only upon the Shannon number $N$, we make use of the
approximation~(\ref{smallK}) for the kernel $D(\rhat,\rhat')$,
and introduce the independent and dependent variable transformations
\begin{equation}
\label{hugesphere}
{\bf x}=\sqrt{\fracd{4\pi}{A}}\,\rhat,\quad
{\bf x}'=\sqrt{\fracd{4\pi}{A}}\,\rhat',\qquad
\psi({\bf x})=g(\rhat),\quad
\psi({\bf x}')=g(\rhat').
\end{equation}
The scaled coordinates ${\bf x},{\bf x}'$ are the
projections of the points $\rhat,\rhat'\in\Omega$ onto a large
sphere $\Omega_{*}$ of squared radius $\|{\bf x}\|^2=4\pi/A$.
The geodesic distance
between the scaled points ${\bf x},{\bf x}'\in\Omega_{*}$
and the differential surface area on $\Omega_{*}$ are
\begin{equation}
\|{\bf x}-{\bf x}'\|=\sqrt{\fracd{4\pi}{A}}\,\Delta\quad
\mbox{and}\quad d\Omega_{*}=\frac{4\pi}{A}\,d\Omega.
\label{Deltatilde}
\end{equation}
Upon making the substitutions~(\ref{hugesphere})--(\ref{Deltatilde}),
equations~(\ref{firsttimeint}) and~(\ref{smallK}) reduce to
\begin{equation}
\label{scaledup}
\int_{R_{*}}\!\!D_{*}({\bf x},{\bf x}')\,\psi({\bf x})\,d\Omega_{*}'
=\lambda\hspace{0.1em}\psi({\bf x}),
\end{equation}
where $R_{*}$ is the projection of the region of concentration $R$
onto the sphere $\Omega_{*}$, and
\begin{equation}
\label{scaledup2}
D_{*}({\bf x},{\bf x}')=\fracd{\sqrt{N}}{2\pi}\,
\fracd{J_1\big(\raisebox{-0.2ex}{$\sqrt{N}$}\,\|{\bf x}-{\bf x}'\|\big)}
{\|{\bf x}-{\bf x}'\|}
\end{equation}
is the symmetric, $N$-dependent Fredholm kernel.

\begin{figure}[b]\centering 
\rotatebox{0}{
\includegraphics[width=0.85\textwidth]{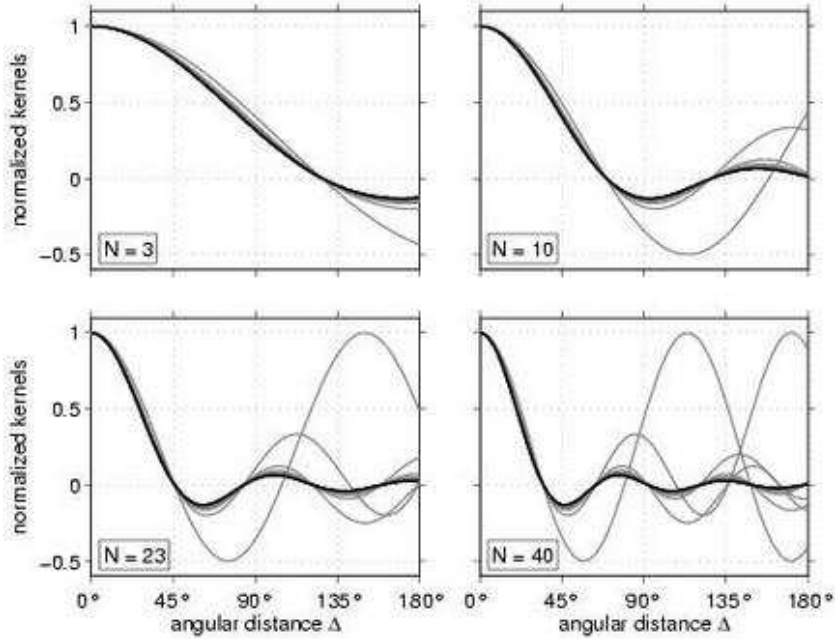}
} 
\caption{\small Comparison of the exact scaled kernels~(\ref{exactscaled})
with the flat-earth asymptotic approximation~(\ref{asympscaled}) (black).
The Shannon number $N=3,10,23,40$ is kept constant in
each of the four panels, and the bandwidth used to compute the
exact scaled kernels varies between
$L=1$ (worst fitting) and $L=100$ (best fitting).}
\label{sdwscaling} 
\end{figure}

Equations~(\ref{scaledup})--(\ref{scaledup2}) are the spherical
analogue of the one-dimensional scaled eigenvalue
equation~(\ref{slepian5}).  The asymptotic eigenvalues
$\lambda_1,\lambda_2,\ldots$ and associated scaled eigenfunctions
$\psi_1({\bf x}),\psi_2({\bf x}),\ldots$ depend upon the maximal
degree $L$ and the area $A$ only through the Shannon number
$N=\Lpot A/(4\pi)$. As in the case of equations~(\ref{firsttimeint})
and~(\ref{eigen2}), we are free to
solve~(\ref{scaledup})--(\ref{scaledup2}) either on all of
$\Omega_{*}$, in which case the eigenfunctions $\psi_1({\bf
x}),\psi_2({\bf x}),\ldots$ are bandlimited, or only in the region of
concentration $R_{*}$, in which case they are spacelimited.  It is
readily verified that the scaling has no effect upon the sum of the
eigenvalues, inasmuch as
\begin{equation}
N=\int_{R_{*}}\!\!D_{*}(0)\,d\Omega_{*}
=\fracd{N}{4\pi}\int_{R_{*}}d\Omega_{*}=N.
\label{shannonstar}
\end{equation}
We show in Appendix~\ref{sec:flatearth} that the scaled eigenvalue
problem~(\ref{scaledup})--(\ref{scaledup2}) is identical to that
governing the two-dimensional concentration problem on a plane. 

The above considerations show that in the limit~(\ref{asylimit}),
we expect the exact Fredholm kernel~(\ref{banddelta}), evaluated
on $\Omega_{*}$ and normalized by its value at zero offset,
\begin{equation}
\fracd{D\big(\sqrt{4\pi/A}\,\Delta\big)}{D(0)}=
\frac{1}{\Lpot}\suml_{l=0}^{L}(2l+1)\,
P_l\!\left(\cos\sqrt{\frac{4\pi}{A}}\,\Delta\right),
\label{exactscaled}
\end{equation}
to be well approximated by the similarly normalized asymptotic kernel
\begin{equation}
\fracd{D_{*}(\Delta)}{D_{*}(0)}=
\frac{2J_1(\sqrt{N}\Delta)}{\sqrt{N}\Delta}.
\label{asympscaled}
\end{equation}
The quality of this asymptotic approximation to the kernel
and the associated flat-earth scaling are illustrated
in Figure~\ref{sdwscaling}. In the four examples shown,
with Shannon numbers $N=3,10,23,40$, the approximation
is excellent even for angular distances as large as $\Delta\approx 135^{\circ}$,
once the maximal spherical harmonic degree exceeds $L=3$\,--\,$4$.

\subsection{Scaled eigenvalue equation for an axisymmetric polar cap}
The flat-earth asymptotic version of the fixed-order colatitudinal eigenvalue
problem~(\ref{Fredholm}) can be obtained in two different ways:
either by an application of the Hilb approximation~(\ref{Hilb})
and the Poisson sum formula~(\ref{Poisson}) to the kernel
$D(\theta,\theta')$ given in equation~(\ref{dlx}), or by using
the addition theorem for Bessel functions \cite[]{Jeffreys+88},
\begin{equation}
J_0(k\Delta)=J_0(k\theta)J_0(k\theta')
+2\suml_{m=1}^{\infty}J_m(k\theta)J_m(k\theta')\cos m(\phi-\phi'),
\label{additionBJ}
\end{equation}
the representation~(\ref{polarg2}) of $g(\theta,\phi)$,
and the orthonormality of the longitudinal functions
$\ldots,\sqrt{2}\cos m\phi,\ldots, 1,\ldots,
\sqrt{2}\sin m\phi,\ldots$ over the interval $0\le\phi\le 2\pi$
to decompose equations~(\ref{firsttimeint}) and~(\ref{smallK})
into a series of individual eigenvalue problems, one for each
order $0\le m\le L$. Using either method, we find that 
equation~(\ref{Fredholm}) can be approximated in the limit~(\ref{asylimit}) by
\begin{equation}
\label{asyfixedm}
\int_0^{\Theta}D(\theta,\theta')\,g(\theta')\,\theta'\,d\theta'
=\lambda\hspace{0.1em}g(\theta),
\end{equation}
where
\begin{equation}
\label{asyfixedm2}
D(\theta,\theta')=\Lpot\int_{0}^{1}\!
J_m\big[(L+1)\hspace{0.1em}p\hspace{0.1em}\theta\big]\,
J_m\big[(L+1)\hspace{0.1em}p\hspace{0.1em}\theta'\big]\,p\,dp.
\end{equation}
It is convenient in the present instance to approximate the area of
the small polar cap by $A=2\pi(1-\cos\Theta)\approx\pi\Theta^2$, and
to introduce scaled coordinates that are slightly different from those
in equations~(\ref{hugesphere}), namely
\begin{equation}
\label{hugesphere2}
x=\theta\hspace{-0.05em}/\Theta,\quad
x'=\theta'\hspace{-0.2em}/\Theta,\qquad
\psi(x)=g(\theta),\quad
\psi(x')=g(\theta').
\end{equation}
This leads to a scaled, fixed-order eigenvalue problem,
\begin{equation}
\label{fixedmpsi}
\int_{0}^{1}D_{*}(x,x')\,\psi(x')\,x'\hspace{0.1em}dx'
=\lambda\hspace{0.1em}\psi(x),
\end{equation}
with an associated kernel
\begin{equation}
\label{fixedmpsi2}
D_{*}(x,x')=4N\int_{0}^{1}\!
J_m\big(2\sqrt{N}\,p\hspace{0.1em}x\big)\,
J_m\big(2\sqrt{N}\,p\hspace{0.1em}x'\big)\,p\,dp,
\end{equation}
whose eigenvalues
$\lambda_1,\lambda_2,\ldots$ and associated scaled
eigenfunctions $\psi_1(x),\psi_2(x),\ldots$ depend
upon the maximal spherical harmonic degree $L$ and the cap radius
$\Theta$ only through the small polar-cap Shannon number
$N=\frac{1}{4}\Lpot\Theta^2$.

Although the polar-cap scaling
relations~(\ref{fixedmpsi})--(\ref{fixedmpsi2})
are strictly valid only in the asymptotic limit
$L\rightarrow\infty,\Theta\rightarrow 0$, the approximation is excellent
even for moderate bandwidths $L$ and sizable cap radii $\Theta$.
For a fixed Shannon number $N=40$, a maximal
degree in the range $25\leq L\leq 40$,
and therefore a cap radius $\Theta=2\sqrt{N}/(L+1)$
in the range $29^{\circ}\geq\Theta\geq 18^{\circ}$,
the agreement between the fixed-order, scaled
eigenfunctions is always within a few percent.
In principle, the asymptotic results~(\ref{fixedmpsi})--(\ref{fixedmpsi2})
would enable the determination of approximate polar cap eigenfunctions
$g(\theta)$ for varying values of $L$ and $\Theta$ by scaling a
pre-computed catalogue of fixed-$N$ eigenfunctions.
In practice, the construction and diagonalization of
the tridiagonal Gr\"{u}nbaum matrix~(\ref{tridiag})
is so straightforward and efficient that it is
preferable to simply compute the optimally concentrated
eigenfunctions $g(\theta)$ exactly.

\subsection{Asymptotic fixed-order Shannon number}
\begin{figure}[h]\centering 
\rotatebox{0}{
\includegraphics[width=0.85\textwidth]{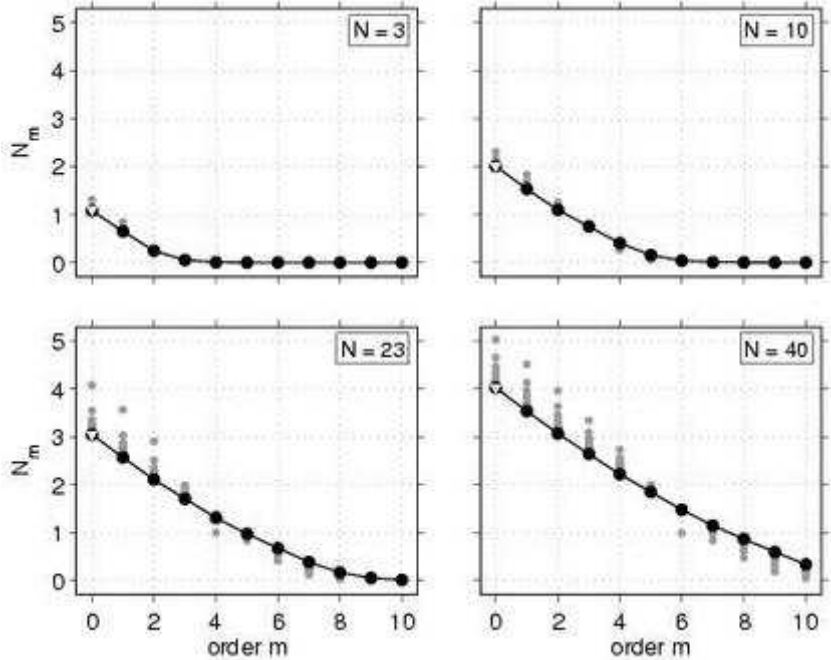}
} 
\caption{\small Comparison of the number $N_m$ of significant
eigenvalues of fixed order $m$ (grey) with the asymptotic
approximation~(\ref{nsubm3}) (black).  The Shannon number
$N=3,10,23,40$ is kept constant in each of the four panels, and the
bandwidth used to compute the exact values of $N_m$ varies between
$L=1$ (worst fitting) and $L=100$ (best fitting). Points inconsistent
with the constraint $A/(4\pi)=N/\Lpot <1$ are not plotted. White
triangles show the simplified zonal ($m=0$) approximation $N_0\approx
(L+1)\Theta/\pi$.}
\label{sdwnsubm} 
\end{figure}

The asymptotic approximation to the number of significant
eigenvalues associated with a given order $m$ is
\ber
N_m&=&\int_0^1D_{*}(x,x)\,x\,dx\nnr\\&=&4N\int_0^1\!\!\int_0^1\!
J_m^2\big(2\sqrt{N}\,p\hspace{0.1em}x\big)\,p\,dp\,x\,dx\nnr\\
&=&{}+2N\left[J^2_m\big(2\sqrt{N}\,\big)+J^2_{m+1}\big(2\sqrt{N}\,\big)\right]\nnr\\
&&{}-(2m+1)\sqrt{N}J_m\big(2\sqrt{N}\,\big)J_{m+1}\big(2\sqrt{N}\,\big)\nnr\\
&&{}-\frac{m}{2}\left[
1-J_0^2\big(2\sqrt{N}\,\big)-2\suml_{n=1}^{m}J_n^2\big(2\sqrt{N}\,\big)
\right].
\label{nsubm3}
\eer
The relationship~(\ref{Nequalsum}) between the total
number $N$ of significant eigenvalues and the number $N_m$ associated with
each order $m$ is preserved in this asymptotic approximation, inasmuch
as 
\ber
N&=& N_0+2\suml_{m=1}^{\infty}N_m\nnr\\
&=&4N\int_{0}^{1}\!\!\int_{0}^{1} \left[J_0^2\big(2\sqrt{N}pq\big)+
2\suml_{m=1}^{\infty}J^2_m\big(2\sqrt{N}pq\big) \right]\,p\,dp\,x\,dx\nnr\\
&=&4N\int_{0}^{1}\!\!\int_{0}^{1}\,p\,dp\,x\,dx=N.  
\eer 
In Figure~\ref{sdwnsubm} we compare the exact fixed-order Shannon numbers
$N_m$, computed by Gauss-Legendre numerical integration of
equation~(\ref{nsubm2}), with the asymptotic result~(\ref{nsubm3}),
for the same values of $N=3,10,23,40$ and $1\leq L\leq 100$
as in Figure~\ref{sdwscaling}.  
The number of significant $m=0$
eigenvalues can be even more simply approximated by
$N_0\approx 2\sqrt{N}/\pi\approx (L+1)\Theta/\pi$, as shown.
\section{Conclusion}
An orthogonal family of bandlimited spherical harmonic expansions that
are optimally concentrated within a finite region $R$ of the unit
sphere can be computed by solving either a symmetric matrix eigenvalue
problem in the spectral domain or an equivalent Fredholm integral
eigenvalue problem in the spatial domain.  Every eigenvalue
$0<\lambda<1$ is a measure of both the spatial concentration of the
bandlimited eigenfunction $g(\rhat)$ and the spectral concentration of
the spacelimited eigenfunction $h(\rhat)$ that coincides with
$g(\rhat)$ inside the region of concentration. The number of
well concentrated eigenfunctions is $N=\Lpot A/(4\pi)$, where $L$ is
the maximal spherical degree and $A$ is the area of the region of
concentration.  Roughly speaking, this Shannon number $N$ is the
dimension of the space of functions $f(\rhat)$ that can be
concentrated both within a finite region $R$ of the sphere and within
a spectral interval $0\le l\le L$. For a small region $A\ll 4\pi$, and
a moderate maximal spherical harmonic degree $L$, the optimally
concentrated bandlimited eigenfunctions $g(\rhat)$ and associated
spacelimited eigenfunctions $h(\rhat)$ can be computed accurately,
even for an irregularly shaped region $R$. In the special, but
important, case of a circular polar cap, every eigenfunction can be
computed accurately, by numerical diagonalization of a commuting
tridiagonal matrix, which has a simple Sturm-Liouville spectrum.  Just
as Slepian's one-dimensional, prolate spheroidal eigentapers have
proven to be extremely useful in time-frequency spectral analysis, we
expect the two-dimensional, spherical eigenfunctions developed here to
have a wide variety of spatiospectral data analysis applications in
fields such as geophysics, planetary science and cosmology.

\section*{Acknowledgments}
F.~J.~S. thanks Ingrid Daubechies, Peter E. Harris, Jean Steiner, Partha
Mitra, Bill Symes, and David Thomson for insightful discussions,
and the D{\'e}partement de G{\'e}ophysique Spatiale et
Plan{\'e}taire at the Institut de Physique du Globe de Paris
for their hospitality. Financial support for this work was
provided by the U.~S.~National Science Foundation under
Grant EAR-0105387. This is IPGP contribution number XXXX. 

\bibliography{bib.bib}
\bibliographystyle{siam}

\appendix
\section{Computational considerations}
\label{sec:Computational}
Here we present a brief description of the numerical methods employed in
this study. All of the computations described here have been performed
using double precision arithmetic. Statements regarding machine
precision refer to double precision, with a roundoff error of
$\sim\! 10^{-16}$.

\subsection{Concentration within a polar cap}
\label{subsec:Comp:polar}
We compute the colatitudinal eigenfunctions
$g_1(\theta),g_2(\theta),\ldots,g_{L-m+1}(\theta)$ of an axisymmmetric
polar cap $0\le\theta\le\Theta$ using three different methods. The
first method is by numerical diagonalization of the
$(L-m+1)\times(L-m+1)$ matrix ${\sfD}$ in equation~(\ref{fixedmeqn}).
We do not implement the Wigner 3-$j$ expression~(\ref{kernel3j}) for
the matrix elements $D_{ll'}$, but instead use Gauss-Legendre
quadrature \cite[]{Press+92} to evaluate the defining
integral~(\ref{kernel4}): \ber D_{ll'}&=&\int_{\cos\Theta}^1
X_{lm}(\arccos\mu)X_{l'm}(\arccos\mu)\,d\mu\nnr\\
&\approx&\sum_{j=1}^Jw_jX_{lm}(\arccos\mu_j) X_{l'm}(\arccos\mu_j),
\label{firstgauss}
\eer
where $\mu_1,\mu_2,\ldots,\mu_J$ are the roots of the Legendre
polynomial $P_J(\bar{\mu})$, rescaled from $-1\le\bar{\mu}_j\le 1$ to
$\cos\Theta\le\mu_j\le 1$, and
\mbox{$w_j=2(1-\bar{\mu}_j^2)^{-1}[P_J^{\prime}(\bar{\mu}_j)]^{-2},j=1,2,\ldots,J$}
are the associated integration weights. Only the uppermost triangular
matrix elements $D_{ll'},l\le l'$ are computed explicitly; the
lowermost elements are infilled using the symmetry $D_{ll'}=D_{l'l}$.
The order $J$ of the Gauss-Legendre integration is adjusted upward
until the $L-m+1$ spatial-domain eigenfunctions
$g_1(\theta),g_2(\theta),\ldots,g_{L-m+1}(\theta)$ satisfy the
orthogonality relations~(\ref{fixedmortho}) to within machine
precision.  The same high-order Gauss-Legendre quadrature rule is used
to evaluate the orthogonality integrals.  The Legendre functions
$X_{lm}(\theta)$ are computed with high accuracy to very high degree
($l\approx 500$) using a recursive algorithm
\cite[]{Libbrecht85,Masters+98}. 

The second method is by numerical solution of the Fredholm
equation~(\ref{Fredholmmu}). Using a Gauss-Legendre
quadrature rule to discretize this equation, we obtain
\begin{equation}
\sum_{j=1}^{J}w_jD(\mu_j,\mu^{\prime}_j)h(\mu^{\prime}_j)=
\lambda\hspace{0.1em}h(\mu_j),\quad j=1,2,\ldots,J.
\label{interpol}
\end{equation}
Equation~(\ref{interpol}) can be rewritten as a symmetric
algebraic eigenvalue equation,
\begin{equation}
(\sfW^{\scriptscriptstyle{1\hspace{-0.01em}/\hspace{-0.01em}2}}
\tilde{\sfD}\sfW^{\scriptscriptstyle{1\hspace{-0.01em}/\hspace{-0.01em}2}})
(\sfW^{\scriptscriptstyle{1\hspace{-0.01em}/\hspace{-0.01em}2}}\tilde{\sfh})=
\lambda\hspace{0.1em}
(\sfW^{\scriptscriptstyle{1\hspace{-0.01em}/\hspace{-0.01em}2}}\tilde{\sfh}),
\label{wonehalf}
\end{equation}
where $\tilde{\sfh}$ is a $J$-dimensional column vector with entries
$\tilde{h}_j=h(\mu_j)$, and where $\tilde{\sfD}$ and $\sfW$ denote the
$J\times J$ matrices with elements
$\tilde{D}_{jj'}=D(\mu_j,\mu^{\prime}_j)$ and
$W_{jj'}=w_j\delta_{jj'}$.  The eigenvalues $\lambda$ and transformed
eigenvectors
$\sfW^{\scriptscriptstyle{1\hspace{-0.01em}/\hspace{-0.01em}2}}\tilde{\sfh}$
are computed by numerical diagonalization of the matrix
$\sfW^{\scriptscriptstyle{1\hspace{-0.01em}/\hspace{-0.01em}2}}
\tilde{\sfD}\sfW^{\scriptscriptstyle{1\hspace{-0.01em}/\hspace{-0.01em}2}}$.
The order of integration $J$ is again chosen to ensure accurate
orthogonality of the spatial-domain eigenfunctions
$h_1(\theta),h_2(\theta),\ldots,h_{L-m+1}(\theta)$.  In the zonal
($m=0$) case the choice $J=L+1$ renders both of the
integrations~(\ref{firstgauss}) and~(\ref{interpol}) exact; for
$m\not= 0$ we use a conservative, larger integration order $J$, since
the integrands are no longer polynomials.  The fixed-order power
spectra shown in Figure~\ref{sdwspectral} are computed by transforming
the spatial-domain eigenfunctions of equation~(\ref{eigen2}) to the
spectral domain, using a spherical harmonic degree range $m\leq l\leq
127$ sufficient to avoid aliasing \cite[]{Sneeuw94}.

Even for moderate values of the maximal degree $L$ and cap radius
$\Theta$, the smallest eigenvalues
$\ldots,\lambda_{L-m},\lambda_{L-m+1}$ fall below machine precision.
The associated, least well concentrated eigenfunctions computed using
either of the above two direct methods are in that case essentially
arbitrary orthogonal members of a numerically degenerate eigenspace,
and are no longer accurate \cite[]{Albertella+99}.  Because of this,
it is not possible to find the optimally excluded eigenfunctions of a
small polar cap, or equivalently the optimally concentrated
eigenfunctions of a large cap, by diagonalization of either of the
matrices $\sfD$ or $\tilde\sfD$.  Fortunately, this difficulty can be
overcome by the third method, which is numerical diagonalization of
the tridiagonal Gr\"{u}nbaum matrix~(\ref{tridiag}). The roughly
equant spacing of the Sturm-Liouville eigenvalues
$\chi_1,\chi_2,\ldots,\chi_{L-m+1}$ enables all of the associated
eigenfunctions to be calculated to within machine precision.  The
spatiospectral concentration factors
$\lambda_1,\lambda_2,\ldots,\lambda_{L-m+1}$ are computed to the same
precision, either by {\it a posteriori\/} matrix multiplication,
$\lambda=\sfg^{\sf{\sst{T}}}{\sfD}\sfg$, or by Gauss-Legendre
integration of the equivalent spatial relation~(\ref{fixedmortho}).
Both the significant and the insignificant eigenvalues computed using
each of the above methods agree to within machine precision, providing
a useful numerical check. Diagonalization of the tridiagonal matrix
$\sfG$ is the only numerically stable way to solve the concentration
problem for either a large polar cap or a large maximal degree $L$. By
an extension of the above analysis, it is even possible to use the
Gr\"{u}nbaum operator ${\mathcal G}$ to compute spacelimited
eigenfunctions $h_{L-m+2}(\theta),h_{L-m+3}(\theta),\ldots$ that are
in the null space \cite[]{Miranian2004}.

\subsection{Concentration within an arbitrarily shaped region}
\label{subsec:Comp:arbit}
We solve the spatiospectral concentration problem for an arbitrarily shaped
region $R$ by numerical diagonalization of the $\Lpot\times\Lpot$
matrix $\sfD$, with elements $D_{lm,l'm'}$
defined by equation~(\ref{Dlmlmpdef}).
Given the splined boundary of $R$,
we first find the northernmost and southernmost points,
with colatitudes $\theta_{\rm n}$ and $\theta_{\rm s}$.
For every $\theta_{\rm n}\le\theta\le\theta_{\rm s}$,
we then find the easternmost and westernmost points,
with longitudes $\phi_{\rm e}(\theta)$ and $\phi_{\rm w}(\theta)$.
In the case of a non-convex region with indentations and protuberances,
there may be several such easternmost and westernmost points,
which we shall index with an additional subscript $i=1,2,\ldots,I$.
The integral over longitude,
\begin{equation}
\Phi_{mm'}(\theta)=\sum_{i=1}^I\int_{\phi_{{\rm w}i}}
^{\phi_{{\rm e}i}}\!\left\{\!\!\begin{array}{c} \cos m\phi \\ \sin m\phi\end{array}\!\!\right\}
\left\{\!\!\begin{array}{c} \cos m'\phi \\ \sin m'\phi\end{array}\!\!\right\}d\phi,
\label{easyone}
\end{equation}
can be done analytically, and we use Gauss-Legendre quadrature to
compute the remaining integral over colatitude:
\ber
D_{lm,l'm'}&=&\int_{\mu_{\rm n}}^{\mu_{\rm s}}
X_{lm}(\arccos\mu)X_{l'm'}(\arccos\mu)\Phi_{mm'}(\arccos\mu)\,d\mu\nnr\\
&\approx&\sum_{j=1}^Jw_jX_{lm}(\arccos\mu_j)X_{l'm'}(\arccos\mu_j)
\Phi_{mm'}(\arccos\mu_j).
\label{hardone}
\eer
As in the case of a polar cap, we adjust the order of the integration
$J$ upward until the spatial-domain eigenfunctions
$g_1(\rhat),g_2(\rhat),\ldots,g_{\Lpot}(\rhat)$
satisfy the orthogonality relations~(\ref{orthog})
to within machine precision. There is no analogue
of the Gr\"{u}nbaum operator ${\mathcal G}$ for
an arbitrarily shaped region, so only the eigenfunctions
associated with eigenvalues that are above machine precision
can be computed accurately.
In most practical applications, this is not a limitation,
since we are generally interested only in the computable, well
concentrated eigenfunctions $g_1(\rhat),g_2(\rhat),\ldots,g_N(\rhat)$,
which are associated with the numerically significant eigenvalues
$\lambda_1,\lambda_2,\ldots,\lambda_N$.

\subsection{Concentration within a non-polar circular cap}
\label{subsec:Comp:away}
One of the principal applications of spherical Slepian functions
in geophysics and planetary physics will be to analyze measurements
within a circularly symmetric region centered upon an arbitrary
geographical location  $\theta_0,\phi_0$
\cite[e.g.,][]{Kido+2003,McGovern+2002,Simons+97b,Simons+97a}.
The preferred procedure for determining the required optimally
concentrated eigenfunctions is first to compute the spherical harmonic
coefficients $g_{lm}$ of the eigenfunctions~(\ref{polarg2})
concentrated within a polar cap $0\leq\theta\leq\Theta$,
and then to rotate these to the desired cap location
\cite[]{Blanco+97,Dahlen+98,Edmonds96,Masters+98}. The actual
windowing of the data for further analysis may either be carried out
in the spectral domain \cite[]{Simons+97a}, or, more simply, by
straightforward multiplication after transformation of
the rotated eigenfunctions to the spatial domain.
If one wishes to avoid spherical harmonic rotation,
it is also possible to compute the rotated eigenfunctions
directly, by performing the numerical integration in
equation~(\ref{easyone}) on the analytically prescribed
boundary of a cap of radius $\Theta$ centered at $\theta_0,\phi_0$,
given by
\begin{equation}
\phi_{\rm w,e}(\theta)=\phi_0\mp\Delta\phi(\theta) \quad \mbox{where} \quad
\Delta\phi(\theta)=\frac{\arccos(\cos\Theta-\cos\theta\cos\theta_0)}
{\sin\theta\sin\theta_0}.
\end{equation}
%We have implemented and intercompared both methods of spatiospectral
%concentration within a non-polar circular cap as a numerical check.

\section{Spatiospectral concentration on a plane}
\label{sec:flatearth}
In one of his many papers extending the one-dimensional analysis,
Slepian \cite[]{Slepian64} considered the spatiospectral concentration
problem in a Cartesian space of arbitrary dimension.  We present a
brief review of the two-dimensional Cartesian concentration problem
here, for comparison with the flat-earth asymptotic analysis presented
in Section~\ref{sec:Asymptotics}.

An arbitrary, real-valued, square-integrable function
$f({\bf r})$ on the plane has the two-dimensional
Fourier representation, analogous to the spherical
harmonic representation~(\ref{expansion}),
\begin{equation}
\label{Bfourier}
f({\bf r})=\left(\fracd{1}{2\pi}\right)^{\!2}\!\!\intinf
F({\bf k})e^{i{\bf k}\cdot{\bf r}}\,d^2{\bf k},\qquad
F({\bf k})=\intinf
f({\bf r})e^{-i{\bf k}\cdot{\bf r}}\,d^2{\bf r}.
\end{equation}
The Fourier orthonormality relation analogous to
equation~(\ref{normalization}) is
\begin{equation}
\label{BFortho}
\left(\fracd{1}{2\pi}\right)^{\!2}\!\!\intinf
e^{i{\bf k}\cdot({\bf r}-{\bf r}')}\,d^2{\bf k}
=\delta({\bf r}-{\bf r}')=\fracd{\delta(\|{\bf r}-{\bf r}'\|)}
{2\pi\|{\bf r}-{\bf r}'\|}.
\end{equation}
Parseval's relation stipulates that the power of any function $f({\bf r})$
in the spatial and spectral domains is identical:
\begin{equation}
\label{Btwopowers}
\intinf f^2({\bf r})\,d^2{\bf r}=\left(\fracd{1}{2\pi}\right)^{\!2}\!\!
\intinf |F({\bf k})|^2\,d^2{\bf k}.
\end{equation}
Equation~(\ref{Btwopowers}) is the planar analogue of
the spherical relation $\|{f}\|^2_{\Omega}=\|\sff\|^2_{\infty}$.

We use $g({\bf r})$ to denote a bandlimited function,
\begin{equation}
\label{Bgdefn}
g({\bf r})=\left(\fracd{1}{2\pi}\right)^{\!2}\!\!
\intK G({\bf k})e^{i\mbf{k}\cdot\mbf{r}}\,d^2{\bf k},
\end{equation}
with no power above a maximal wavenumber $K$.
By analogy with the optimization criterion~(\ref{normratio}),
we seek to concentrate the power of $g({\bf r})$ into a finite region
$R$:
\begin{equation}
\label{Bnormratio}
\lambda=\fracd{\int_R g^2\,d^2{\bf r}}
{\intinf g^2\,d^2{\bf r}}=\mbox{maximum}.
\end{equation}
Bandlimited functions $g({\bf r})$ that maximize the
ratio $\lambda$ in equation~(\ref{Bnormratio}) are
solutions to the Fourier domain eigenvalue equation,
analogous to equation~(\ref{eigen1}),
\begin{equation}
\label{Beigen1}
\intK D({\bf k},{\bf k}')\,G({\bf k})\,d^2{\bf k}'
=\lambda\hspace{0.1em}G({\bf k}),\quad\|{\bf k}\|\le K,
\end{equation}
where
\begin{equation}
\label{Beigen2}
D({\bf k},{\bf k}')=\left(\fracd{1}{2\pi}\right)^{\!2}\!\!
\int_R e^{i({\bf k}-{\bf k}')\cdot{\bf r}}\,d^2{\bf r}.
\end{equation}
The corresponding problem in the spatial domain,
analogous to equation~(\ref{firsttimeint}), is
\begin{equation}
\label{Beigen3}
\int_R D({\bf r},{\bf r}')\,g({\bf r})\,d^2{\bf r}'
=\lambda\hspace{0.1em}g({\bf r}),\quad|{\bf r}|\le\infty,
\end{equation}
where
\begin{equation}
\label{Beigen4}
D({\bf r},{\bf r}')=\left(\fracd{1}{2\pi}\right)^{\!2}\!\!
\intK e^{i{\bf k}\cdot({\bf r}-{\bf r}')}\,d^2{\bf k}.
\end{equation}
Spacelimited eigenfunctions $h({\bf r})$, which vanish outside
of the region $R$, satisfy the same eigenvalue equation~(\ref{Beigen3}),
but with the domain of solution properly restricted:
\begin{equation}
\label{Beigen5}
\int_R D({\bf r},{\bf r}')\,h({\bf r})\,d^2{\bf r}'
=\lambda\hspace{0.1em}h({\bf r}),\quad{\bf r}\in R.
\end{equation}
The associated eigenvalue $0<\lambda <1$ is a measure of both the spatial concentration
of $g({\bf r})$ within the region $R$ and the spectral concentration
of $h({\bf r})$ within the wavenumber range $\|{\bf k}\|\le K$.

For consistency with~(\ref{eigorder}), we rank order the eigenvalues
so that $\lambda_1\ge\lambda_2\ge\ldots$.  The bandlimited, spatial-domain
eigenfunctions $g_1({\bf r}),g_1({\bf r}),\ldots$ may be
chosen to be orthonormal over the whole plane $\|{\bf r}\|\le\infty$
and orthogonal over the region $R$:
\begin{equation}
\intinf g_{\alpha}g_{\beta}\,d^2{\bf r}=\delta_{\alpha\beta}\quad\mbox{and}\quad
\int_R g_{\alpha}g_{\beta}\,d^2{\bf r}=\lambda_{\alpha}\delta_{\alpha\beta}.
\label{Borthog}
\end{equation}
Ths sum of the eigenvalues, or Shannon number, is given by
\begin{equation}
\label{Bshannon}
N=\sum_{\alpha}^{\infty}\lambda_{\alpha}=\intK D({\bf k},{\bf k})\,d^2{\bf k}
=\int_R D({\bf r},{\bf r})\,d^2{\bf r}=K^2\,\fracd{A}{4\pi},
\end{equation}
where $A$ is the area of the region of concentration $R$.
Equations~(\ref{Borthog}) and~(\ref{Bshannon}) are the planar
analogues of the spherical relations~(\ref{orthog}) and~(\ref{shannon}).

Comparison of equations~(\ref{BFortho}) and~(\ref{Beigen4})
shows that the spatial-domain kernel $D({\bf r},{\bf r}')$,
like its spherical counterpart $D(\rhat,\rhat')$, is a
bandlimited delta function.
Upon introducing polar coordinates and integrating over the
angle, we can reduce $D({\bf r},{\bf r}')$ to
a form reminiscent of the representation~(\ref{smallK}):
\begin{equation}
\label{Bbesskern}
D({\bf r},{\bf r}')=
\fracd{1}{2\pi}\int_0^K\!J_0(k\|{\bf r}-{\bf r}'\|)\,k\,dk\nnr
=\fracd{K\hspace{0.1em}
J_1(K\|{\bf r}-{\bf r}'\|)}{2\pi \|{\bf r}-{\bf r}'\|}.
\end{equation}
Upon introducing scaled independent and dependent variables
analogous to~(\ref{hugesphere}),
\begin{equation}
\label{Bhuge}
{\bf x}=\sqrt{\fracd{4\pi}{A}}\,{\bf r},\quad
{\bf x}'=\sqrt{\fracd{4\pi}{A}}\,{\bf r}',\qquad
\psi({\bf x})=g({\bf r}),\quad
\psi({\bf x}')=g({\bf r}'),
\end{equation}
we can rewrite equations~(\ref{Beigen3}) and~(\ref{Bbesskern})
in a form identical to~(\ref{scaledup})--(\ref{scaledup2})
and analogous to~(\ref{slepian5}):
\begin{equation}
\frac{K}{2\pi}\int_{R_{*}}\!\fracd{J_1(K\|{\bf x}-{\bf x}'\|)}
{\|{\bf x}-{\bf x}'\|}\,\psi({\bf x})=\lambda\hspace{0.1em}\psi({\bf x}),
\label{Bscaled}
\end{equation}
where $R_{*}$ is the image of the region of concentration $R$
under the mapping~(\ref{Bhuge}).

If the region of concentration $R$ is a circle of radius $Q$,
then a polar coordinate, ${\bf r}=(q,\phi)$, representation
analogous to~(\ref{polarg2}),
\begin{equation}
g(q,\phi)=\left\{
\begin{array}{l@{\quad\mbox{if}\hspace{0.6em}}l}
\rule[-2mm]{0mm}{6mm}\sqrt{2}\,g(q)\cos m\phi & -L\le m<0\\
\rule[-2mm]{0mm}{6mm}g(q)                     & m=0\\
\rule[-2mm]{0mm}{6mm}\sqrt{2}\,g(q)\sin m\phi & 0< m\le L,\\
\end{array}
\right.
\label{Bpolarg}
\end{equation}
may be used to decompose equations~(\ref{Beigen3}) and~(\ref{Bbesskern})
into a series of fixed-order eigenvalue problems analogous
to~(\ref{asyfixedm})--(\ref{asyfixedm2}):
\begin{equation}
\int_{0}^{Q}D(q,q')\,g(q')\,q'\,dq'
={\lambda}\hspace{0.1em}g(q),
\label{Bfixedm}
\end{equation}
where
\begin{equation}
D(q,q')=K^2\!\int_{0}^{1}
J_m(Kp\hspace{0.05em}q)\,
J_m(Kp\hspace{0.05em}q')\,p\,dp.
\label{Bfixedm2}
\end{equation}
The transformations
\begin{equation}
\label{Bmapping2}
x=q/Q,\quad
x'=q'\hspace{-0.2em}/Q,\qquad
\psi(x)=g(q),\quad
\psi(x')=g(q')
\end{equation}
convert equations~(\ref{Bfixedm})--(\ref{Bfixedm2})
into the scaled eigenvalue problem
\begin{equation}
\label{Bscaled2}
4N\int_{0}^{1}\!\!\int_{0}^{1}\!
J_m\big(2\sqrt{N}\,p\hspace{0.05em}x\big)\,
J_m\big(2\sqrt{N}\,p\hspace{0.05em}x'\big)\,p\,dp\,\psi(x')\,x'\hspace{0.1em}dx'
=\lambda\hspace{0.1em}\psi(x),
\end{equation}
which is identical to~(\ref{fixedmpsi})--(\ref{fixedmpsi2}),
and dependent only upon the Shannon number $N=\frac{1}{4}K^2Q^2$.
Slepian \cite[]{Slepian64} has noted that equation~(\ref{Bscaled2})
is an iterated version of the equivalent ``square root'' equation
\begin{equation}
2\sqrt{N}\int_0^1\!J_m\big(2\sqrt{N}\,x\hspace{0.05em}x'\big)\,
\psi(x')\,x'\hspace{0.1em}dx'=\sqrt{\lambda}\,\psi(x).
\label{slepsqrt}
\end{equation}
The eigenvalues $\lambda_1,\lambda_2,\ldots$ and
eigenfunctions $\psi_1(x),\psi_2(x),\ldots$ of
equation~(\ref{Bscaled2}) may alternatively be
found by solving the equivalent equation~(\ref{slepsqrt}).

In the asymptotic limit~(\ref{asylimit}), both the
general and axisymmetric spherical concentration problems
are seen to be identical to the corresponding concentration problem in a plane,
with the maximal wave\-number $K$ replaced by the integer
$L+1$. The planar problem exhibits exact Shannon-number
scaling analogous to that of the one-dimensional problem~(\ref{slepian5}),
whereas the scaling of the spherical problem is only asymptotic.
Equation~(\ref{nsubm3}) giving the number of significant
eigenvalues $N_m$ associated with each angular order $m$
is exact in the case of the plane.

\label{lastpage}

\end{document}